\documentclass[prb,showpacs,preprintnumbers,amsmath,amssymb,twosides]{revtex4-1}

\usepackage[applemac]{inputenc}
\usepackage{amsfonts}
\usepackage{graphicx}% Include figure files
\usepackage{dcolumn}% Align table columns on decimal point
\usepackage{bm}% bold math
\usepackage[english]{babel}

\newcommand{\bea}{\begin{eqnarray}}
\newcommand{\eea}{\end{eqnarray}}

\newcommand{\beqn}{\begin{eqnarray}}
\newcommand{\eeqn}{\end{eqnarray}}

\def\d#1d#2{\frac{\partial #1}{\partial #2}}

\usepackage{amssymb,mathrsfs,euscript}
\usepackage{amssymb}
\usepackage{amsmath}
\usepackage{amscd}
\usepackage{mathrsfs,euscript}
\usepackage[matrix]{xypic}
\usepackage{makeidx}
\makeindex

\textheight9in  \def\DATE{\today}
\newtheorem{theorem}{Theorem}
\newtheorem{definition}[theorem]{Definition}
\newtheorem{corollary}[theorem]{Corollary}

\newtheorem{lemma}{Lemma}

\newtheorem{proposition}{Proposition}

\newcommand\R{\mathbb{R}}

\def\Lie{{\mathcal{L}{\it Lie}}}
\def\pre-Lie{{\mbox {\it pre-\Lie}}}

\def\Gip3ass{{\mbox {$\mathcal{G}_i$\it{-}$p^3$\it{ass}}}}
\def\Gp3ass#1{{\mbox {$\mathcal{G}_#1$\it{-}$p^3$\it{ass}}}}

 \catcode`\@=13

\def\Lie{{\mathcal{L}{\it Lie}}}
\def\pre-Lie{{\mbox {\it pre-\Lie}}}

 \def\Lie{\hbox{{$\mathcal L$}{\it ie\/}}}

\usepackage{graphicx}

\begin{document}

\title{A new algebraic and arithmetic framework for interval computations}
\author{Dr Abdel KENOUFI}
\email{transmat.consulting@gmail.com}
\affiliation{TranSmaT Scientific Consulting, Strasbourg, France}
\author{Dr Nicolas GOZE}
\email{nicolas.goze@uha.fr}
\affiliation{Laboratoire de Math\'ematiques, Informatique et Applications, Universit\'e de Haute de Haute-Alsace, Mulhouse, France}
\author{Dr Elisabeth REMM}
\email{elisabeth.remm@uha.fr}
\affiliation{Laboratoire de Math\'ematiques, Informatique et Applications, Universit\'e de Haute de Haute-Alsace, Mulhouse, France}
\author{Pr Michel GOZE}
\email{michel.goze@uha.fr}
\affiliation{Laboratoire de Math\'ematiques, Informatique et Applications, Universit\'e de Haute de Haute-Alsace, Mulhouse, France}\date{\today}

\begin{abstract}
In this paper we propose some very promissing results in interval arithmetics which permit to build well-defined arithmetics including distributivity of multiplication and division according addition and substraction. Thus, it allows to build all algebraic operations and functions on intervals. This will avoid completely the wrapping effects and data dependance. Some simple applications for matrix eigenvalues calculations, inversion of symmetric matrices and finally optimization are exhibited in the object-oriented programming language \emph{python}.
\end{abstract}

\maketitle
\section{History}
The first mathematician who has used intervals was the famous Archimedes from Syracuse (287-212 b.C). He has proposed a two-sides bounding of $\pi$ : $3+\frac{10}{71}<\pi<3+\frac{1}{7}$ using polygons and a systematic method to improve it. In the beginning of the twentieth century, the American
mathematician and physicist Wiener, published two papers \cite{Wiener1,Wiener2}, and used intervals to give an interpretation to the position and the time of a system. More papers on the subject were written \cite{Dwyer,Sunaga,Warmus1,Warmus2} only after Second World War. Nowadays, we consider R.E. Moore \cite{Moore1,Moore2,Moore3,Moore4,Moore5} as the first mathematician who has proposed a framework for interval arithmetics and analysis.
The interval arithmetic, or interval analysis has been introduced to compute very quickly range bounds (for example if a data is given up to an incertitude). Now interval arithmetic is a computing system which permits to perform error analysis by computing mathematic bounds. The extensions of the areas of applications are important: non linear problems, PDE, inverse problems. It finds a large place of
applications in controllability, automatism, robotics, embedded systems, biomedical, haptic interfaces, form optimization, analysis of architecture plans, ...\\
Interval calculations are used nowadays as a powerful tool for global optimization and set inversion\cite{Jaulin}. Several groups have developped some softwares and libraries to perform those new apparoaches such as \emph{INTLAB}\cite{intlab}, \emph{INTOPT90} and \emph{GLOBSOL}\cite{intopt90-globsol},\emph{Numerica}\cite{numerica}. But our goal for this article, is not to replace the semantic approach of intervals, which has to be adapted to each problem by the engineer or the scientist, but to propose a new arithmetic of intervals, which allows to avoid the wrapping and data dependance effects. It yields to a better construction of inclusion functions.\\
We expose in this paper the main results of a PhD thesis\cite{nicolas} defended by one the author and some consecutive numerical applications. The plan is the following. In a first time we define a real Banach structure on the completion $\overline{\mathbb{IR}}$ of the semi group of intervals $\mathbb{IR}$, with a vector space structure. This permits to  define the notion of differential function with values of $%
\overline{\mathbb{IR}}$ and to use
some important tools and the fixed point theorem. Next we extend the classical product to have a
distributivity property.   With this approach we obtain a notion of differential calculus and a natural linear algebra on the set of intervals.
After that, we gives some examples in a \emph{python} implementation and we end this article by giving some simple numerical applications : optimization of interval functions, interval matrix diagonalization, and inversion of symmetric matrices .

\section{An algebraic approach to the set of intervals.}

In this section we present the set of intervals as a normed vector space. We
define also a four-dimensional associative algebra whose product gives the
product of intervals in any cases.

\subsection{Minkowski operations}
An interval is a bounded non empty connected closed subset of $\mathbb{R}.$ Let $\mathbb{IR}$ be the set of intervals.
The semantical arithmetic operations on intervals, called Minkowski operations,  are defined such that the
result of the corresponding operation on elements belonging to operand
intervals belongs to the resulting interval. That is, if $\diamond $ denotes
one of the semantical operations $+,-,\ast $, we have,
if $X$ and $\ Y$ are bounded intervals of $\mathbb{R}$,
\begin{equation*}
X\diamond Y=\{x\diamond y\text{ }/\text{ }x\in X,\text{ }y\in Y\},
\end{equation*}

 In many problems using interval arithmetic, that is the set $\mathbb{IR}$ with the Minkowski operations, there exists an
informal transfers principle which permits, to associate with a real
function $f$ a function define on the set of intervals $\mathbb{IR}$ which
coincides with $f$ on the interval reduced to a point. But this transferred
function is not unique. For example, if we consider the real function $%
f(x)=x^{2}+x=x(x+1),$ we associate naturally the functions $\widetilde{f}%
_{1}:\mathbb{IR}\longrightarrow \mathbb{IR}$ given by $\widetilde{f}%
_{1}(X)=X(X+1)$ and $\widetilde{f}_{2}(X)=X^{2}+X.$ These two functions do not coincide. Usually this
problem is removed considering the most interesting transfers. But the
qualitative "interesting" depends of the studied model and it is not given
by a formal process. In this section, we determine a natural extension $\overline{%
\mathbb{IR}}$ of $\mathbb{IR}$ provided with a vector space structure. The
vectorial substraction $X\smallsetminus Y$ does not correspond to the
semantical difference of intervals and the interval  $\smallsetminus X$ has no real
interpretation. But these "negative" intervals have a computational role. If
a problem conduce to a "negative" result, then this problem is "pervert"
(see Lazare Carnot with his feeling on the natural negative number).

Let $\mathbb{IR}$ be the set of intervals. It is in one to one
correspondence with the half plane of $\mathbb{R}^{2}$:
\begin{equation*}
\mathcal{P}_{1}=\{(a,b),a\leq b\}.
\end{equation*}%
This set is closed for the addition and $\mathcal{P}_{1}$ is endowed with a
regular semi-group structure. Let $\mathcal{P}_{2}$ be the half plane
symmetric to $\mathcal{P}_{1}$ with respect to the first bisector $\Delta $
of equation $y-x=0.$ The substraction on $\mathbb{IR}$, which is not the
symmetric operation of $+$, corresponds to the following operation on $%
\mathcal{P}_{1}$:
\begin{equation*}
(a,b)-(c,d)=(a,b)+s_{\Delta }\circ s_{0}(c,d),
\end{equation*}%
where $s_{0}$ is the symmetry with respect to $0$, and $s_{\Delta }$ with
respect to $\Delta .$ The multiplication $\ast $ is not globally defined.
Consider the following subset of $\mathcal{P}_{1}$:
\begin{equation*}
\left\{
\begin{array}{l}
\mathcal{P}_{1,1}=\{(a,b)\in \mathcal{P}_{1},a\geq 0,b\geq 0\}, \\
\mathcal{P}_{1,2}=\{(a,b)\in \mathcal{P}_{1},a\leq 0,b\geq 0\}, \\
\mathcal{P}_{1,3}=\{(a,b)\in \mathcal{P}_{1},a\leq 0,b\leq 0\}. \\
\end{array}%
\right.
\end{equation*}%
We have the following cases:

1) If $(a,b),(c,d)\in \mathcal{P}_{1,1}$ the product is written $(a,b)\ast
(c,d)=(ac,bd).$

\noindent

\noindent The vectors $\ e_{1}=(1,1)$ and $e_{2}=(0,1)$ generate $\mathcal{P}%
_{1,1}$ that is any $(x,y)$ in $\mathcal{P}_{1,1}$, can be decomposed as
\begin{equation*}
(x,y)=xe_{1}+(y-x)e_{2},\text{with }x>0\text{ and }y-x>0.
\end{equation*}%
The multiplication corresponds in this case to the following associative
commutative algebra:
\begin{equation*}
\left\{
\begin{array}{l}
e_{1}e_{1}=e_{1}, \\
e_{1}e_{2}=e_{2}e_{1}=e_{2}e_{2}=e_{2}.%
\end{array}%
\right.
\end{equation*}

2) Assume that $(a,b)\in \mathcal{P}_{1,1}$ and $(c,d)\in \mathcal{P}_{1,2}$
so $c\leq 0$ and $d\geq 0.$ Thus we obtain $(a,b)\ast (c,d)=(bc,bd)$ and
this product does not depend of $a.$ Then we obtain the same result for any $%
a<b$. The product $(a,b)\ast (c,d)=(bc,bd)$ corresponds to
\begin{equation*}
\left\{
\begin{array}{l}
e_{1}e_{1}=e_2e_1=e_1 \\
e_{1}e_{2}=e_2e_2=e_{2}
\end{array}
\right.
\end{equation*}
This algebra is not commutative and it is different from the previous.

3) If $(a,b)\in \mathcal{P}_{1,1}$ and $(c,d)\in \mathcal{P}_{1,3}$ then $%
a\geq 0,b\geq 0$ and $c\leq 0,d\leq0$ and we have $(a,b)\ast (c,d)=(bc,ad).$
Let $e_{1}=(1,1)$, $e_{2}=(0,1)$. This product corresponds to the following
associative algebra:
\begin{equation*}
\left\{
\begin{array}{l}
e_{1}e_{1}=e_{1}, \\
e_{1}e_{2}=e_{2}, \\
e_{2}e_{1}=e_{1}-e_{2}.
\end{array}
\right.
\end{equation*}
This algebra is not associative because $(e_2e_1)e_1 \neq e_2(e_1e_1)$.
We have similar results for the cases $(\mathcal{P}_{1,2},\mathcal{P}
_{1,2}),\ (\mathcal{P}_{1,2},\mathcal{P}_{1,3})$ and $(\mathcal{P}_{1,3},
\mathcal{P}_{1,3}).$

An objective of this paper is to present an associative algebra which contains all these results.

\subsection{The real vector space $\overline{\mathbb{IR}}$}

We recall briefly the construction proposed by Markov \cite{Ma} to define a
structure of abelian group. As $(\mathbb{IR},+)$ is a commutative and
regular semi-group, the quotient set, denoted by $\overline{\mathbb{IR}}$,
associated with the equivalence relations:
\begin{equation*}
(x,y)\sim (z,t)\Longleftrightarrow x+t=y+z,
\end{equation*}%
for all $x,y,z,t\in \mathbb{IR}$, is provided with a structure of abelian
group for the natural addition:
\begin{equation*}
\overline{(x,y)}+\overline{(z,t)}=\overline{(x+z,y+t)}
\end{equation*}%
where $\overline{(x,y)}$ is the equivalence class of $(x,y)$. We denote by $%
\smallsetminus \overline{(x,y)}$ the opposite of $\overline{(x,y)}$ . We
have $\smallsetminus \overline{(x,y)}=$ $\overline{(y,x)}$. If $x=[a,a]$, $%
a\in \mathbb{R}$, then $\overline{(x,0)}=\overline{(0,-x)}$ where $%
-x=[-a,-a] $, and $\smallsetminus \overline{(x,0)}=\overline{(0,x)}$. In
this case, we identify $x=[a,a]$ with $a$ and we denote always by $\mathbb{R}
$ the subset of intervals of type $[a,a]$. Naturally, the group $\overline{\mathbb{IR}}$
is isomorphic to the additive group $\mathbb{R}^2$ by the isomorphism $(%
\overline{([a,b],[c,d])}\rightarrow (a-c,b-d)$. We find the notion of
generalized interval.
\begin{proposition}
Let $\mathcal{X}=\overline{(x,y)}$ be in $\overline{\mathbb{IR}}$. Thus
\begin{enumerate}
\item If $l(y)<l(x),$ there is an unique $A\in \mathbb{IR\setminus R}$ such
that $\mathcal{X}=\overline{(A,0)},$
\item If $l(y)>l(x),$ there is an unique $A\in \mathbb{IR}\setminus \mathbb{%
R }$ such that $\mathcal{X}=\overline{(0,A)}=\smallsetminus \overline{(A,0)}%
, $
\item If $l(y)=l(x),$ there is an unique $A=\alpha \in \mathbb{R}$ such that
$\mathcal{X}=\overline{(\alpha ,0)}=\overline{(0,-\alpha )}. $
\end{enumerate}
\end{proposition}
Any element $\mathcal{X}=\overline{(A,0)}$ with $A\in \mathbb{IR-R}$ is said
positive and we write $\mathcal{X}>0.$ Any element $\mathcal{X}=\overline{%
(0,A)}$ with $A\in \mathbb{IR-R}$ is said negative and we write $\mathcal{X}%
<0.$ We write $\mathcal{X}\geq \mathcal{X}^{\prime }$ if $\mathcal{X}%
\smallsetminus \mathcal{X}^{\prime }\geq 0.$ For example if $\mathcal{X}$ \ and $\mathcal{X}^{\prime }$ are
positive, $
\mathcal{X}\geq \mathcal{X}^{\prime }\Longleftrightarrow l(\mathcal{X})\geq
l(\mathcal{X}^{\prime }).
$. The elements $\overline{(\alpha ,0)}$ with $\alpha \in \mathbb{R}
^{\ast }$ are neither positive nor negative.

\medskip

In \cite{Ma}, one defines on the abelian group $\overline{\mathbb{IR}}$ , a
structure of quasi linear space.
Our approach is a little bit different. We propose to construct a real
vector space structure.  We consider the
external multiplication:
\begin{equation*}
\cdot :\mathbb{R}\times \overline{\mathbb{IR}}\longrightarrow \overline{%
\mathbb{IR}}
\end{equation*}%
defined, for all $A\in \mathbb{IR}$, by
\begin{equation*}
\left\{
\begin{array}{c}
\alpha \cdot \overline{(A,0)}\text{\ }=\overline{(\alpha A,0)}\text{ ,} \\
\alpha \cdot \overline{(0,A)}\text{\ }=\overline{(0,\alpha A)}\text{ ,}%
\end{array}%
\right.
\end{equation*}%
for all $\alpha >0.$ If $\alpha <0$ we put $\beta =-\alpha $. So we put:
\begin{equation*}
\left\{
\begin{array}{c}
\alpha \cdot \overline{(A,0)}\text{\ }=\overline{(0,\beta A)}, \\
\alpha \cdot \overline{(0,A)}\text{\ }=\overline{(\beta A,0)}.%
\end{array}%
\right.
\end{equation*}

\noindent We denote $\alpha \mathcal{X}$ instead of $\alpha \cdot \mathcal{X}
.$ This operation satisfies

\begin{enumerate}
\item For any $\alpha \in \mathbb{R}$ and $\mathcal{X}\in \overline{\mathbb{%
IR}}$ we have:
\begin{equation*}
\left\{
\begin{array}{c}
\alpha (\smallsetminus \mathcal{X})=\smallsetminus (\alpha \mathcal{X}), \\
(-\alpha )\mathcal{X}=\smallsetminus (\alpha \mathcal{X}).%
\end{array}
\right.
\end{equation*}

\item For all $\alpha ,\beta \in \mathbb{R}$, and for all $\mathcal{X},%
\mathcal{X} ^{\prime }\in \overline{\mathbb{IR}}$, we have
\begin{equation*}
\left\{
\begin{tabular}{l}
$(\alpha +\beta )\mathcal{X}=\alpha \mathcal{X}+\beta \mathcal{X}$, \\
$\alpha (\mathcal{X}+\mathcal{X}^{\prime })=\alpha \mathcal{X}+\alpha
\mathcal{X}^{\prime }$, \\
$(\alpha \beta )\mathcal{X}=\alpha (\beta \mathcal{X}).$%
\end{tabular}
\right.
\end{equation*}
\end{enumerate}

\begin{theorem}
The triplet $(\overline{\mathbb{IR}},+,\cdot )$ is a real vector space and
the vectors  $\mathcal{X}_{1}=\overline{([0,1],0)}$ and $\mathcal{\ X}_{2}=%
\overline{([1,1],0)}$ of $\overline{\mathbb{IR}}$ determine a basis of $%
\overline{\mathbb{IR}}.$ So $\dim _{\mathbb{R}}\overline{\mathbb{IR}}=2.$
\end{theorem}

\noindent \textit{Proof. }We have the following decompositions:%
\begin{equation*}
\left\{
\begin{array}{c}
\overline{([a,b],0)\text{ }}=(b-a)\mathcal{X}_{1}+a\mathcal{X}_{2}, \\
\overline{(0,[c,d])\text{ }}=(c-d)\mathcal{X}_{1}-c\mathcal{X}_{2}.%
\end{array}
\right.
\end{equation*}
The linear map
\begin{equation*}
\varphi :\overline{\mathbb{IR}}\longrightarrow \mathbb{R}^{2}
\end{equation*}
defined by
\begin{equation*}
\left\{
\begin{array}{l}
\varphi (\,\overline{([a,b],0)}\text{ })=(b-a,a), \\
\varphi (\,\overline{(0,[c,d])}\text{ })=(c-d,-c)%
\end{array}
\right.
\end{equation*}
is a linear isomorphism and $\overline{\mathbb{IR}}$ is canonically
isomorphic to $\mathbb{R}^{2}$.

\medskip

\noindent \textbf{Remark. }Let $E$ be the subspace generated by $\mathcal{X}
_{2}.$ The vectors of $E$ correspond to the elements which have a non
defined sign. Then the relation $\leq $ defined in the paragraph $1.2$ gives
an order relation on the quotient space $\overline{\mathbb{IR}}/E.$

\subsection{A Banach structure on $\overline{\mathbb{IR}}$}

 Any element $%
\mathcal{X}\in \overline{\mathbb{IR}}$ $\ $is written $\overline{(A,0)}$ or $%
\overline{(0,A)}.$\ We define its length $l(\mathcal{X})$ as the length of $%
A $ and its center as $c(A)$ or $-c(A)$ in the second case.

\begin{theorem}
The map $||$ $||$ $:$ $\overline{\mathbb{IR}}\longrightarrow $ $\mathbb{R}$
given by
\begin{equation*}
||\mathcal{X}||=l(\mathcal{X})+|c(\mathcal{X})|
\end{equation*}
for any $\mathcal{X}\in $ $\overline{\mathbb{IR}}$ is a norm$.$
\end{theorem}

\noindent \noindent \noindent \textit{Proof}. We have to verify the
following axioms:
\begin{equation*}
\left\{
\begin{array}{l}
1)\text{ }||\mathcal{X}||=0\Longleftrightarrow \mathcal{X}=0, \\
2)\text{ }\forall \lambda \in \mathbb{R}\text{ }||\lambda \mathcal{X}
||=|\lambda |||\mathcal{X}||,\text{ } \\
3)\text{ }||\mathcal{X}+\mathcal{X}^{\prime }||\leq ||\mathcal{X}||+||
\mathcal{X}^{\prime }||.%
\end{array}
\right.
\end{equation*}

\noindent 1) If $||\mathcal{X}||=0$, then $l(\mathcal{X})=|c(\mathcal{X})|=0$
and $\mathcal{X}=0.$

\bigskip

\noindent \noindent \noindent \noindent 2) Let $\lambda \in \mathbb{R}.$ We
have
\begin{equation*}
||\lambda \mathcal{X}||=l(\lambda \mathcal{X})+|c(\lambda \mathcal{X}
)|=|\lambda |l(\mathcal{X})+|\lambda ||c(\mathcal{X})|=|\lambda |||\mathcal{%
X }||.
\end{equation*}

\medskip

\noindent \noindent 3) We consider that $I$ refers to $\mathcal{X}$ and $J$
refers to $\mathcal{X}^{\prime}$ thus $\mathcal{X}=\overline{(I,0)}$ or $=
\overline{(0,I)}$. We have to study the two different cases:

\noindent i) If $\mathcal{X}+\mathcal{X}^{\prime }=\overline{(I+J,0)}$ or $%
\overline{(0,I+J)}$, then
\begin{eqnarray*}
||\mathcal{X}+\mathcal{X}^{\prime }||
&=&l(I+J)+|c(I+J)|=l(I)+l(J)+|c(I)+c(J)|\leq l(I)+|c(I)|+l(J)+|c(J)| \\
&=&||\mathcal{X}||+||\mathcal{X}^{\prime }||.
\end{eqnarray*}

\noindent ii) Let $\mathcal{X}+\mathcal{X}^{\prime }=\overline{(I,J)}.$ If $%
\overline{(I,J)}=\overline{(K,0)}$ then $K+J=I$ and
\begin{equation*}
||\mathcal{X}+\mathcal{X}^{\prime }||=||\overline{(K,0)}
||=l(K)+|c(K)|=l(I)-l(J)+|c(I)-c(J)|
\end{equation*}
that is
\begin{equation*}
||\mathcal{X}+\mathcal{X}^{\prime }||\leq l(I)+|c(I)|-l(J)+|c(J)|\leq
l(I)+|c(I)|+l(J)+|c(J)|=||\mathcal{X}||+||\mathcal{X}^{\prime }||.
\end{equation*}
So we have a norm on $\overline{\mathbb{IR}}.$

\begin{theorem}
The normed vector space $\overline{\mathbb{IR}}$ is a Banach space.
\end{theorem}

\noindent \textit{Proof}. In fact, all the norms on $\mathbb{R}^2$ are
equivalent and $\mathbb{R}^2$ is a Banach space for any norm. The vector
space $\overline{\mathbb{IR}}$ is isomorphic to $\mathbb{R}^2$. Thus it is
complete.

\medskip

\noindent\textbf{Remarks. }

\begin{enumerate}
\item To define the topology of the normed space $\overline{\mathbb{IR}}$,
it is sufficient to describe the $\mathcal{\varepsilon }$-neighborhood of
any point $\chi _{0}\in $ $\overline{\mathbb{IR}}$ for $\varepsilon $ a
positive infinitesimal number. We can give a geometrical representation,
considering $\chi _{0}=\overline{([a,b],0)}$ represented by the point $%
(a,b)\in \mathbb{R}^{2}.$ We assume that $\chi _{0}=\overline{([a,b],0)}$
and $\mathcal{\varepsilon }$ an infinitesimal real number. Let $A_{1},\cdots
,A_{4}$ the points $A_{1}=(a-\varepsilon ,b-\varepsilon ),A_{2}=(a+\frac{%
\varepsilon }{2},b-\frac{\varepsilon }{2}),A_{3}=(a+\varepsilon
,b+\varepsilon ),A_{4}=(a-\frac{\varepsilon }{2},b+\frac{\varepsilon }{2})$.
If $0<a<b$, then the $\mathcal{\varepsilon }$-neighborhood of $\chi _{0}=%
\overline{([a,b],0)}$ is represented by the parallelograms whose vertices
are $A_{1},A_{2},A_{3},A_{4}.$

\item We can consider another equivalent norms on $\overline{\mathbb{IR}}$.
For example
\begin{equation*}
||\mathcal{X}||=||\smallsetminus \mathcal{X}||=Sup(|x|,|y|)
\end{equation*}%
where $\mathcal{X}=\overline{([x,y],0)}$. But we prefer the initial one
because it has a better geometrical interpretation.
\end{enumerate}

\section{\protect\bigskip Differential calculus on $\overline{\mathbb{IR}}
$}

As $\overline{\mathbb{IR}}$ is a Banach space, we can describe a notion of
differential function on it. \ Consider $\mathcal{X} _{0}=\overline{%
(X_{0},0) }$ in $\overline{\mathbb{IR}}$ . The norm $||.||$ defines a
topology on $\overline{\mathbb{IR}}$ whose a basis of neighborhoods is given
by the balls $\mathcal{B}(X_{0},\varepsilon )=\{X\in \overline{\mathbb{IR}}%
,|| \mathcal{X} \smallsetminus \mathcal{X} _{0}||<\varepsilon \}.$ Let us
characterize the elements of $\mathcal{B}(X_{0},\varepsilon ).$ $\mathcal{X}
_{0}=\overline{(X_{0},0)}=\overline{([a,b],0)}.$

\begin{proposition}
Consider $\mathcal{X}_{0}=\overline{(X_{0},0)}=\overline{([a,b],0)}$ in $%
\overline{\mathbb{IR}}$ and $\varepsilon \simeq 0$, $\varepsilon >0$. Then
every element of $\mathcal{B}(X_{0},\varepsilon )$ is of type $\mathcal{X}=%
\overline{(X,0)}$ and satisfies
\begin{equation*}
l(X)\in B_{\mathbb{R}}(l(X_{0}),\varepsilon _{1})\text{ and }c(X)\in B_{%
\mathbb{R}}(c(X_{0}),\varepsilon _{2})
\end{equation*}%
with $\varepsilon _{1},\varepsilon _{2}\geq 0$ and $\varepsilon _{1}+%
\mathcal{\varepsilon }_{2}\leq \varepsilon ,$ where $B_{\mathbb{R}}(x,a)$ is
the canonical open ball in $\mathbb{R}$\ of center $x$ and radius $a.$
\end{proposition}

\noindent\textit{Proof. }\textit{First case : }Assume that $\mathcal{X} =
\overline{(X,0)}=\overline{([x,y],0)}$ . We have
\begin{eqnarray*}
\mathcal{X} \smallsetminus \mathcal{X} _{0} &=&\overline{(X,X_{0})}=
\overline{([x,y],[a,b])} \\
&=&\left\{
\begin{array}{c}
\overline{([x-a,y-b],0)}\text{ if }l(X)\geq l(X_{0}) \\
\overline{(0,[a-x,b-y])}\text{ if }l(X)\leq l(X_{0})%
\end{array}
\right.
\end{eqnarray*}
If $l(X)\geq l(X_{0})$ we have
\begin{eqnarray*}
||\mathcal{X} \smallsetminus \mathcal{X} _{0}|| &=&(y-b)-(x-a)+\left\vert
\frac{y-b+x-a}{2}\right\vert \\
&=&l(X)-l(X_{0})+|c(X)-c(X_{0})|.
\end{eqnarray*}
As \ $l(X)-l(X_{0})\geq 0$ and $|c(X)-c(X_{0})|\geq 0,$ each one of this
term if less than $\varepsilon.$ If $l(X)\leq l(X_{0})$ we have
\begin{equation*}
||\mathcal{X} \smallsetminus \mathcal{X}
_{0}||=l(X_{0})-l(X)+|c(X_{0})-c(X)|.
\end{equation*}
and we have the same result.

\bigskip

\noindent\textit{Second case : }Consider $\mathcal{X} =\overline{(0,X)}=
\overline{([x,y],0)}$ .\ We have
\begin{equation*}
\mathcal{X} \smallsetminus \mathcal{X} _{0}=\overline{(0,X_{0}+X)}=\overline{
([x+a,y+b])}
\end{equation*}
and
\begin{equation*}
||\mathcal{X} \smallsetminus \mathcal{X}
_{0}||=l(X_{0})+l(X)+|c(X_{0})+c(X)|.
\end{equation*}
In this case, we cannot have $||\mathcal{X} \smallsetminus \mathcal{X}
_{0}||<\varepsilon$ thus $X\notin \mathcal{B}(X_{0},\varepsilon).$

\medskip

\begin{definition}
A function $f:\overline{\mathbb{IR}}\longrightarrow \overline{\mathbb{IR}}$
is continuous at $\mathcal{X}_{0}$ if
\begin{equation*}
\forall \varepsilon >0,\exists \eta >0\text{ such that }||\mathcal{X}%
\smallsetminus \mathcal{X}_{0}||<\eta \text{ implies }||f(\mathcal{X}%
)\smallsetminus f(\mathcal{X}_{0})||<\varepsilon .
\end{equation*}
\end{definition}

Consider $(\mathcal{X} _{1},\mathcal{X} _{2})$ the basis of $\overline{
\mathbb{IR}}$ given in section 2. We have
\begin{equation*}
f(\mathcal{X} )=f_{1}(\mathcal{X} )\mathcal{X} _{1}+f_{2}(\mathcal{X} )
\mathcal{X} _{2}\text{ with }f_{i}:\overline{\mathbb{IR}}\longrightarrow
\mathbb{R}\text{.}
\end{equation*}
If $f$ is continuous at $\mathcal{X} _{0}$ so
\begin{equation*}
f(\mathcal{X} )\smallsetminus f(\mathcal{X} _{0})=(f_{1}(\mathcal{X}
)-f_{1}( \mathcal{X} _{0}))\mathcal{X} _{1}+(f_{2}(\mathcal{X} )-f_{2}(%
\mathcal{X} _{0}))\mathcal{X} _{2}.
\end{equation*}
To simplify notations let $\alpha =$ $f_{1}(\mathcal{X} )-f_{1}(\mathcal{X}
_{0})$ and $\beta $ =$f_{2}(\mathcal{X} )-f_{2}(\mathcal{X} _{0}).$ If $||f(
\mathcal{X} )\smallsetminus f(\mathcal{X} _{0})||<\varepsilon$, and if we
assume $f_{1}(\mathcal{X} )-f_{1}(\mathcal{X} _{0})>0$ and $f_{2}(\mathcal{X}
)-f_{2}(\mathcal{X} _{0})>0$ (other cases are similar), then we have
\begin{equation*}
l(\alpha \mathcal{X} _{1}+\beta \mathcal{X} _{2})=l\overline{([\beta ,\alpha
+\beta ],0)}<\varepsilon
\end{equation*}
thus $f_{1}(\mathcal{X} )-f_{1}(\mathcal{X} _{0})<\varepsilon.$ Similarly,
\begin{equation*}
c(\alpha \mathcal{X} _{1}+\beta \mathcal{X} _{2})=c\overline{([\beta ,\alpha
+\beta ],0)}=\frac{\alpha }{2}+\beta <\varepsilon
\end{equation*}
and this implies that $f_{2}(\mathcal{X} )-f_{2}(\mathcal{X}
_{0})<\varepsilon.$

\begin{corollary}
$f$ is continuous at $\mathcal{X} _{0}$ if and only if $f_{1}$ and $f_{2}$
are continuous at $\mathcal{X} _{0}.\medskip $
\end{corollary}

\begin{definition}
Consider $\mathcal{X}_{0}$ in $\overline{\mathbb{IR}}$ and $f:\overline{
\mathbb{IR}}\longrightarrow \overline{\mathbb{IR}}$ continuous. We say that $%
f$ is differentiable at $\mathcal{X}_{0}$ if there is $g:\overline{\mathbb{%
IR }}\longrightarrow \overline{\mathbb{IR}}$ linear such as
\begin{equation*}
||f(\mathcal{X})\smallsetminus f(\mathcal{X}_{0})\smallsetminus g(\mathcal{X}
\smallsetminus \mathcal{X}_{0})||=o(||\mathcal{X}\smallsetminus \mathcal{X}
_{0}||).
\end{equation*}
\end{definition}

\noindent \textbf{Examples. }

\begin{itemize}
\item $f(\mathcal{X} )=\mathcal{X}$. This function is continuous at any point and differentiable. It's derivative is $f'(\mathcal{X} )=1$.
\item $f(\mathcal{X} )=\mathcal{X} ^{2}.$ Consider $\mathcal{X} _{0}=
\overline{(X_{0},0)}=\overline{([a,b],0)}$ and $\mathcal{X} \in \mathcal{B}
(X_{0},\varepsilon).$ We have
\begin{eqnarray*}
||\mathcal{X} ^{2}\smallsetminus \mathcal{X} _{0}^{2}|| &=&||(\mathcal{X}
\smallsetminus \mathcal{X} _{0})(\mathcal{X} +\mathcal{X} _{0})|| \\
&\leq &||\mathcal{X} \smallsetminus \mathcal{X} _{0}||||\mathcal{X} +
\mathcal{X} _{0}||.
\end{eqnarray*}
Given $\varepsilon>0,$ let $\eta =\dfrac{\varepsilon}{||\mathcal{X} +
\mathcal{X} _{0}||},$ thus if $||\mathcal{X} \smallsetminus \mathcal{X}
_{0}||<\eta $, we have $||\mathcal{X} ^{2}\smallsetminus \mathcal{X}
_{0}^{2}||<\varepsilon$ and $f$ is continuous and differentiable. It is easy to prove that $f'(\mathcal{X} )=2\mathcal{X}$ is its derivative.

\item Consider $P=a_{0}+a_{1}X+\cdots +a_{n}X^{n}\in \mathbb{R[X]}$. We
define $f:\overline{\mathbb{IR}}\longrightarrow \overline{\mathbb{IR}}$ with
$f(\mathcal{X} )=a_{0}\mathcal{X} _{2}+a_{1}\mathcal{X} +\cdots +a_{n}^{n}
\mathcal{X} ^{n}$ where $\mathcal{X} ^{n}=\mathcal{X} \cdot\mathcal{X}
^{n-1} $ . From the previous example, all monomials are continuous and differentiable, it implies that $f$ is continuous and differentiable as well.

\item Consider the function $Q_2$ given by $Q_2([x,y])=[x^2, y^2]$ if $|x|<|y|$ and
$Q_2([x,y]=[y^2,x^2]$  in the other case. This function is not differentiable.
\end{itemize}

\section{A 4-dimensional associative algebra associated with $
\overline{\mathbb{IR}}$}

%\subsection{Definition of $\mathcal{A}_{4}$}

In introduction, we have observed that the semi-group $\mathbb{IR}$ is
identified to $\mathcal{P}_{1,1}\cup \mathcal{P}_{1,2}\cup \mathcal{P}_{1,3}.
$ Let us consider the following vectors of $\mathbb{R}^{2}$

\begin{equation*}
\left\{
\begin{array}{l}
e_{1}=(1,1), \\
e_{2}=(0,1), \\
e_{3}=(-1,0), \\
e_{4}=(-1,-1).%
\end{array}%
\right.
\end{equation*}%
They correspond to the intervals $[1,1],[0,1],[-1,0],[-1,-1].$ Any point of $%
\mathcal{P}_{1,1}\cup \mathcal{P}_{1,2}\cup \mathcal{P}_{1,3}$ admits the
decomposition

\begin{equation*}
(a,b)=\alpha _{1}e_{1}+\alpha _{2}e_{2}+\alpha _{3}e_{3}+\alpha _{4}e_{4}
\end{equation*}%
with $\alpha _{i}\geq 0.$The dependance relations between the vectors $e_{i}$
are%
\begin{equation*}
\left\{
\begin{array}{l}
e_{2}=e_{3}+e_{1} \\
e_{4}=-e_{1}.%
\end{array}%
\right.
\end{equation*}%
Thus there exists a unique decomposition of $(a,b)$ in a chosen basis such
that the coefficients are non negative. These basis are $\{e_{1,}e_{2}\}$
for $\mathcal{P}_{1,1},$ $\{e_{2},e_{3}\}$ for $\mathcal{P}_{1,2},$ $%
\{e_{3},e_{4}\}$ for $\mathcal{P}_{1,3},$ \ Let us consider the free algebra
of basis $\{e_{1},e_{2},e_{3},e_{4}\}$ whose products correspond to the
Minkowski products. The multiplication table is%
\begin{equation*}
\begin{tabular}{|l|l|l|l|l|}
\hline
& $e_{1}$ & $e_{2}$ & $e_{3}$ & $e_{4}$ \\ \hline
$e_{1}$ & $e_{1}$ & $e_{2}$ & $e_{3}$ & $e_{4}$ \\ \hline
$e_{2}$ & $e_{2}$ & $e_{2}$ & $e_{3}$ & $e_{3}$ \\ \hline
$e_{3}$ & $e_{3}$ & $e_{3}$ & $e_{2}$ & $e_{2}$ \\ \hline
$e_{4}$ & $e_{4}$ & $e_{3}$ & $e_{2}$ & $e_{1}$ \\ \hline
\end{tabular}%
\ .
\end{equation*}%
This algebra is associative. Let $\varphi : \overline{\mathbb{IR}} \rightarrow \mathcal{A}_{4}$ the natural injective embedding. If we identify an interval with its image in $\mathcal{A}_{4}$, we have:

\begin{theorem}
The multiplication of intervals in the algebra $\mathcal{A}_{4}$ is
distributive with respect the addition.
\end{theorem}

 The application  is not bijective. Its image on the elements $\mathcal{X}=
\overline{(x,0)}=\overline{([a,b],0)}$  is:
\begin{equation*}
\left\{
\begin{array}{l}
x=[a,b]\in \mathcal{P}_{1,1},\varphi (\mathcal{X})=ae_{1}+(b-a)e_{2}\quad (a\geq
0,b-a\geq 0) \\
x=[a,b]\in \mathcal{P}_{1,2},\varphi (\mathcal{X})=-ae_{3}+be_{2}\quad (-a\geq 0,b\geq
0) \\
x=[a,b]\in \mathcal{P}_{1,3},\varphi (\mathcal{X})=-be_{4}+(b-a)e_{3}\quad (-b\geq
0,b-a\geq 0).%
\end{array}%
\right.
\end{equation*}%
Consider in $\mathcal{A}_{4}$ the linear subspace $F$ generated by the vectors $%
e_{1}-e_{2}+e_{3},e_{1}+e_{4}.$ As%
\begin{equation*}
\begin{array}{l}
(e_{1}+e_{4})(e_{1}+e_{4})=2(e_{1}+e_{4}) \\
(e_{1}+e_{4})(e_{1}-e_{2}+e_{3})=e_{1}+e_{4} \\
(e_{1}-e_{2}+e_{3})(e_{1}-e_{2}+e_{3})=e_{1,}%
\end{array}%
\end{equation*}%
$F$ is not a subalgebra of $\mathcal{A}_{4}.$ Let us consider the map%
\begin{equation*}
\overline{\varphi }:\mathbb{IR\rightarrow }\mathcal{A}_{4}/F
\end{equation*}%
defined from $\varphi $ and the canonical projection on the quotient vector
space $\mathcal{A}_{4}/F$. A vector $x=\sum \alpha _{i}e_{i}\in \mathcal{A}%
_{4}$ is equivalent to a vector of $\mathcal{A}_{4}$ with positive
components if and only if
\begin{equation*}
\alpha _{2}+\alpha _{3}\geq 0.
\end{equation*}%
In this case, all the vectors equivalent to $x=\sum \alpha _{i}e_{i}$ with $%
\alpha _{2}+\alpha _{3}\geq 0$ correspond to the interval $[\alpha
_{1}-\alpha _{3}-\alpha _{4},\alpha _{1}+\alpha _{2}-\alpha _{4}]$ of $%
\mathbb{IR}$. Thus we have for any equivalent classes of $\mathcal{A}_{4}/F$
associated with $\sum \alpha _{i}e_{i}$ with $\alpha _{2}+\alpha _{3}\geq 0$
a preimage in $\mathbb{IR}$. The map $\overline{\varphi }$ is
injective. In fact, two intervals belonging to pieces $\mathcal{P}_{1,i},%
\mathcal{P}_{1,j}$ with $i\neq j$, have distinguish images. Now if $(a,b)$
and $(c,d)$ belong to the same piece, for example $\mathcal{P}_{1,1}$, thus%
\begin{equation*}
\overline{\varphi }(a,b)=\left\{ (a+\lambda +\mu ,b-a-\lambda ,\lambda ,\mu
),\lambda ,\mu \in \mathbb{R}\text{.}\right\}
\end{equation*}%
If $\overline{\varphi }(c,d)=\overline{\varphi }(a,b)$, there are $\lambda
,\mu \in \mathbb{R}$ such that $(c,d)=(a+\lambda +\mu ,b-a-\lambda ,\lambda
,\mu )$. This gives $a=c,b=d$. We have the same results for all the other
pieces.Thus $\overline{\varphi }:\mathbb{IR\rightarrow }\mathcal{A}_{4}/F$
is bijective on its image, that is the hyperplane of \ $\mathcal{A}_{4}/F$
corresponding to $\alpha _{2}+\alpha _{3}\geq 0$.

Practically the multiplication  of two intervals will so be made: let $%
X,Y\in \mathbb{R}$. Thus $X=\sum \alpha _{i}e_{i},Y=\sum \beta _{i}e_{i}$
with $\alpha _{i},\beta _{j}\geq 0$ and we have the product
\begin{equation*}
X\bullet Y=\overline{\varphi }^{-1}(\varphi (X).\varphi (Y))
\end{equation*}%
this product is well defined because $\overline{\varphi (X).\varphi (Y)}\in
Im\overline{\varphi }.$\ This product is distributive because%
\begin{equation*}
\begin{array}{cl}
X\bullet (Y+Z) & =\overline{\varphi }^{-1}(\varphi (X).\varphi (Y+Z)) \\
& =\overline{\varphi }^{-1}(\varphi (X).(\varphi (Y)+\varphi (Z)) \\
& =\overline{\varphi }^{-1}(\varphi (X).\varphi (Y)+\varphi (X).\varphi (Z))
\\
& =X\bullet Y+X\bullet Z%
\end{array}%
\end{equation*}

\bigskip

\noindent \textbf{Remark.} We have%
\begin{equation*}
\overline{\varphi }^{-1}(\varphi (X).\varphi (Y+Z))\neq \overline{\varphi }%
^{-1}(\varphi (X)).\overline{\varphi }^{-1}(\varphi (Y+Z))).
\end{equation*}%
We shall be careful not to return in $\mathbb{IR}$ during the calculations
as long as the result is not found. Otherwise we find the semantic problems
of the distributivity.

\bigskip

We extend naturally the map $\varphi :\mathbb{IR\rightarrow }\mathcal{A}_{4}$
to $\overline{\mathbb{IR}}$ by%
\begin{equation*}
\left\{
\begin{array}{l}
\varphi \overline{(A,0)}=\varphi (A) \\
\varphi \overline{(0,A)}=-\varphi (A)%
\end{array}%
\right.
\end{equation*}%
for every $A\in \mathbb{IR}$.

\begin{theorem}
The multiplication%
\begin{equation*}
\mathcal{X}^{\prime }\bullet \mathcal{X}^{\prime \prime }=\overline{\varphi }%
^{-1}(\varphi (\mathcal{X}^{\prime }).\varphi (\mathcal{X}^{\prime \prime }))
\end{equation*}%
is distributive with respect the addition.
\end{theorem}

\noindent\textit{Proof. }This is a direct consequence of the previous computations.

%\subsection{Algebraic study of $\mathcal{A}_{4}$}

\bigskip

In $\mathcal{A}_{4}$ we consider the change of basis
\begin{equation*}
\left\{
\begin{array}{l}
e_{1}^{\prime }=e_{1}-e_{2} \\
e_{i}^{\prime }=e_{i},i=2,3 \\
e_{4}^{\prime }=e_{4}-e_{3}.%
\end{array}%
\right.
\end{equation*}%
This change of basis shows that $\mathcal{A}_{4}$ is isomorphic to $\mathcal{%
A}_{4}^{\prime }$
\begin{equation*}
\begin{tabular}{|l|l|l|l|l|}
\hline
& $e_{1}$ & $e_{2}$ & $e_{3}$ & $e_{4}$ \\ \hline
$e_{1}$ & $e_{1}$ & $0$ & $0$ & $e_{4}$ \\ \hline
$e_{2}$ & $0$ & $e_{2}$ & $e_{3}$ & $0$ \\ \hline
$e_{3}$ & $0$ & $e_{3}$ & $e_{2}$ & $0$ \\ \hline
$e_{4}$ & $e_{4}$ & $0$ & $0$ & $e_{1}$ \\ \hline
\end{tabular}%
\ .
\end{equation*}%
The unit of $\mathcal{A}_{4}^{\prime }$ is the vector $e_{1}+e_{2}.$ This
algebra is a direct sum of two ideals: $\mathcal{A}_{4}^{\prime }=I_{1}+I_{2}
$ where $I_{1}$ is generated by $e_{1}$ and $e_{4}$ and $I_{2}$ is generated
by $e_{2}$ and $e_{3}.$ It is not an integral domain, that is, we have
divisors of $0.$ For example $e_{1}\cdot e_{2}=0.$

\begin{proposition}
The multiplicative group $\mathcal{A}_{4}^{\ast }$ \ of invertible elements of $\mathcal{A}_4$ is the set of elements $%
x=(x_{1},x_{2},x_{3},x_{4})$ such that
\begin{equation*}
\left\{
\begin{array}{c}
x_{4}\neq \pm x_{1}, \\
x_{3}\neq \pm x_{2}.%
\end{array}
\right.
\end{equation*}
If $x\in $ $\mathcal{A}_{4}^{\ast }$ we have:
\begin{equation*}
x^{-1}=\left( \frac{x_{1}}{x_{1}^{2}-x_{4}^{2}},\frac{x_{2}}{
x_{2}^{2}-x_{3}^{2}},\frac{x_{3}}{x_{2}^{2}-x_{3}^{2}},\frac{x_{4}}{
x_{1}^{2}-x_{4}^{2}}\right) .
\end{equation*}
\end{proposition}

%\subsection{Monotony property}

\bigskip

Let us  compute the product of intervals using the product in $%
\mathcal{A}_{4}$ and we compare with the Minkowski product. Let $X=[a,b]$
and $Y=[c,d]$ two intervals.

\begin{lemma}
If $X$ and $Y$ are not in the same piece $\mathcal{P}_{1,i}$, then $X\bullet Y$
corresponds to the Minkowski product.
\end{lemma}

\noindent\textit{Proof. \ }i) If\textit{\ }$X\in \mathcal{P}_{1,1}$ and $Y\in
\mathcal{P}_{1,2\text{ }}$then $\varphi (X)=(a,b-a,0,0)$ and $\varphi
(Y)=(0,d,-c,0).$\ Thus%
\begin{equation*}
\begin{array}{cl}
\varphi (X)\varphi (Y) & =(ae_{1}+(b-a)e_{2})(de_{2}-ce_{3}) \\
& =bde_{2}-cbe_{3} \\
& =(0,bd,-cb,0) \\
& =\varphi ([cb,bd]).%
\end{array}%
\end{equation*}%
ii) If $X\in \mathcal{P}_{1,1}$ and $Y\in \mathcal{P}_{1,3\text{ }}$then $%
\varphi (X)=(a,b-a,0,0)$ and $\varphi (Y)=(0,0,d-c,-d).$\ Thus%
\begin{equation*}
\begin{array}{cl}
\varphi (X)\varphi (Y) & =(ae_{1}+(b-a)e_{2})((d-c)e_{3}-de_{4}) \\
& =(ad-bc)e_{3}-ade_{4} \\
& =(0,0,ad-cb,-ad) \\
& =\varphi ([bc,ad]).%
\end{array}%
\end{equation*}%
iii) If $X\in \mathcal{P}_{1,2}$ and $Y\in \mathcal{P}_{1,3\text{ }}$then $%
\varphi (X)=(0,b,-a,0)$ and $\varphi (Y)=(0,0,d-c,-d).$\ Thus%
\begin{equation*}
\begin{array}{cl}
\varphi (X)\varphi (Y) & =(be_{2}-ae_{3})((d-c)e_{3}-de_{4}) \\
& =ace_{2}-bce_{3} \\
& =(0,ac,-cb,0) \\
& =\varphi ([bc,ad]).%
\end{array}%
\end{equation*}

\begin{lemma}
If $X$ an $Y$ are both in the same piece $\mathcal{P}_{1,1}$ or $\mathcal{P}_{1,3}$, then the product $X\bullet Y$ corresponds to the Minkowski product.
\end{lemma}
The proof is analogous to the previous.

\noindent

Let us assume that $X=[a,b]$ and $Y=[c,d]$ belong to  $\mathcal{P}_{1,2}$. Thus $\varphi (X)=(0,b,-a,0)$ and $\varphi (Y)=(0,d,-c,0).$ We obtain
$$XY=(be_2-ae_3)(de_2-ce_3)=(bd+ac)e_2+(-bc-ad)e_3.$$
Thus
$$[a,b][c,d]=[bc+ad,bd+ac].$$
This result is greater that all the possible results associated with the Minkowski product. However, we have the following property:

\begin{proposition}
\textbf{Monotony property: }Let $\mathcal{X}_{1},\mathcal{X}_{2}\in
\overline{\mathbb{IR}}$. Then
\begin{equation*}
\left\{
\begin{array}{l}
\mathcal{X}_{1}\subset \mathcal{X}_{2}\Longrightarrow \mathcal{X}_{1}\bullet
\mathcal{Z}\subset \mathcal{X}_{2}\bullet \mathcal{Z}\text{ for all }
\mathcal{Z}\in \overline{\mathbb{IR}}. \\
\overline{\varphi }(\mathcal{X}_{1})\leq \overline{\varphi }(\mathcal{X}
_{2})\Longrightarrow \overline{\varphi }(\mathcal{X}_{1}\bullet \mathcal{Z}
)\leq \overline{\varphi }(\mathcal{X}_{2}\bullet \mathcal{Z})%
\end{array}
\right.
\end{equation*}
\end{proposition}
The order relation on $\mathcal{A}_{4}$ that ones uses here is
\begin{equation*}
\left\{
\begin{array}{l}
(x_{1},x_{2},0,0)\leq (y_{1},y_{2},0,0)\Longleftrightarrow y_{1}\leq x_{1}
\text{ and }x_{2}\leq y_{2}, \\
(x_{1},x_{2},0,0)\leq (0,y_{2},y_{3},0)\Longleftrightarrow \text{ }x_{2}\leq
y_{2}, \\
(0,x_{2},x_{3},0)\leq (0,y_{2},y_{3},0)\Longleftrightarrow x_{3}\leq y_{3}
\text{ and }x_{2}\leq y_{2}, \\
(0,0,x_{3},x_{4})\leq (0,y_{2},y_{3},0)\Longleftrightarrow \text{ }x_{3}\leq
y_{3}, \\
(0,0,x_{3},x_{4})\leq (0,0,y_{3},y_{4})\Longleftrightarrow x_{3}\leq y_{3}
\text{ and }y_{4}\leq x_{4}.%
\end{array}
\right.
\end{equation*}

\noindent \textit{Proof. }Let us note that the second property is equivalent
to the first.\ It is its translation in $\overline{\mathcal{A}_{4}}.$
We can suppose that $\mathcal{X}_{1}$ and $\mathcal{X}_{2}$ are intervals belonging moreover to $\mathcal{P}_{1,2}$:
$\varphi(\mathcal{X}_{1})=(0,b,-a,0), \varphi(\mathcal{X}_{2})=(0,d,-c,0)$. If $\varphi(\mathcal{Z})=(z_1,z_2,z_3,z_4)$, then
\begin{equation*}
\left\{
\begin{array}{l}
\overline{\varphi }(\mathcal{X}_{1}\bullet \mathcal{Z}
)=(0,bz_1+bz_{2}-az_{3}-az_4,-az_1+bz_{3}-az_{2}+bz_4,0), \\
\overline{\varphi }(\mathcal{X}_{2}\bullet \mathcal{Z}
)=(0,dz_1+dz_{2}-cz_{3}-cz_4,-cz_1+dz_{3}-cz_{2}+dz_4,0).%
\end{array}
\right.
\end{equation*}
Thus
\begin{equation*}
\overline{\varphi }(\mathcal{X}_{1}\bullet \mathcal{Z})\leq \overline{
\varphi }(\mathcal{X}_{2}\bullet \mathcal{Z})\Longleftrightarrow \left\{
\begin{array}{c}
(b-d)(z_1+z_{2})-(a-c)(z_{3}-z_4)\leq 0, \\
-(a-c)(z_1+z_2)+(b-d)(z_{3}=z_4)\leq 0.%
\end{array}
\right.
\end{equation*}
But $(b-d)$, $-(a-c)\leq 0$ and $\ z_{2},z_{3}\geq 0$. This
implies $\overline{\varphi }(\mathcal{X}_{1}\bullet \mathcal{Z})\leq
\overline{\varphi }(\mathcal{X}_{2}\bullet \mathcal{Z}).$

\section{\protect\bigskip The algebras $\mathcal{A}_{n}$ and an better
result of the product}
We can refine our result of the product to come closer to the result of Minkowski. Consider the one dimensional extension $\mathcal{A}_4 \oplus \R{e_5}=\mathcal{A}_5$, where $e_5$ is a vector corresponding to the interval $[-1,1]$ of $\mathcal{P}_{1,2}$. The multiplication table of $\mathcal{A}_5$ is
\begin{equation*}
\begin{tabular}{|l|l|l|l|l|l|}
\hline
& $e_{1}$ & $e_{2}$ & $e_{3}$ & $e_{4}$ & $e_5$\\ \hline
$e_{1}$ & $e_{1}$ & $e_{2}$ & $e_{3}$ & $e_{4}$ & $e_5$ \\ \hline
$e_{2}$ & $e_{2}$ & $e_{2}$ & $e_{3}$ & $e_{3}$ & $e_5$ \\ \hline
$e_{3}$ & $e_{3}$ & $e_{3}$ & $e_{2}$ & $e_{2}$& $e_5$\\ \hline
$e_{4}$ & $e_{4}$ & $e_{3}$ & $e_{2}$ & $e_{1}$ & $e_5$ \\ \hline
$e_{5}$ & $e_{5}$ & $e_{5}$ & $e_{5}$ & $e_{5}$ & $e_5$ \\ \hline
\end{tabular}%
\ .
\end{equation*}%
The piece $\mathcal{P}_{1,2}$ is written $\mathcal{P}_{1,2}=\mathcal{P}_{1,2,1}\cup \mathcal{P}_{1,2,1} $ where $\mathcal{P}_{1,2,1}=\{[a,b], -a\leq b \}$ and $\mathcal{P}_{1,2,2}=\{[a,b], -a\geq  b\}$.
If $X=[a,b] \in \mathcal{P}_{1,2,1}$ and $Y=[c,d] \in \mathcal{P}_{1,2,2}$, thus
$$\varphi(X).\varphi(Y)=(0,b+a,0,0,-a).(0,0,-c-d,0,d)=(0,-(a+b)(c+d),0,0,a(c+d)+bd).$$
Thus we have
$$X\bullet Y=[-bd-ac-ad,-bc].$$

\bigskip

\noindent{\bf Example}
Let $X=[-2,3]$ and $Y=[-4,2]$. We have $X \in \mathcal{P}_{1,2,1}$ and $Y \in \mathcal{P}_{1,2,2}$. The product in $ \mathcal{A}_4$ gives
$$X \bullet Y=[-16,14].$$
The product in $ \mathcal{A}_5$ gives
$$X \bullet Y=[-12,10].$$
The Minkowski product is
$$[-2,3].[-4,2]=[-12,8].$$
Thus the product in $ \mathcal{A}_5$  is better.

\noindent {\bf Conclusion.} Considering a partition of $\mathcal{P}_{1,2}$, we can define an extension of $\mathcal{A}_4$ of dimension $n$, the choice of $n$ depends  on the approach wanted of the Minkowski product. For example, let us consider the vector $e_6$ corresponding to the interval $[-1,\frac{1}{2}]$. Thus the Minkowsky product gives $e_6.e_6=e_7$ where $e_7$ corresponds to $[-\frac{1}{2},1]$. We obtain a $7$-dimensional associative algebra whose table of multiplication is
\begin{equation*}
\begin{tabular}{|l|l|l|l|l|l|l|l|}
\hline
& $e_{1}$ & $e_{2}$ & $e_{3}$ & $e_{4}$ & $e_5$& $e_{6}$ & $e_{7}$ \\ \hline
$e_{1}$ & $e_{1}$ & $e_{2}$ & $e_{3}$ & $e_{4}$ & $e_5$ & $e_{6}$ & $e_{7}$ \\ \hline
$e_{2}$ & $e_{2}$ & $e_{2}$ & $e_{3}$ & $e_{3}$ & $e_5$ & $e_{6}$ & $e_{7}$ \\ \hline
$e_{3}$ & $e_{3}$ & $e_{3}$ & $e_{2}$ & $e_{2}$& $e_5$ & $e_{7}$ & $e_{6}$\\ \hline
$e_{4}$ & $e_{4}$ & $e_{3}$ & $e_{2}$ & $e_{1}$ & $e_5$ & $e_{7}$ & $e_{6}$\\ \hline
$e_{5}$ & $e_{5}$ & $e_{5}$ & $e_{5}$ & $e_{5}$ & $e_5$ & $e_{5}$ & $e_{5}$\\ \hline
$e_6$   & $e_{6}$ & $e_{6}$ & $e_{7}$ & $e_{7}$ & $e_5$ & $e_{7}$ & $e_{6}$ \\ \hline
$e_7$   & $e_{7}$ & $e_{7}$ & $e_{6}$ & $e_{6}$ & $e_5$ & $e_{6}$ & $e_{7}$ \\ \hline
\end{tabular}%
\ .
\end{equation*}%
\noindent{\bf Example}
Let $X=[-2,3]$ and $Y=[-4,2]$. The decomposition on the basis $\{e_1,\cdots,e_7\}$ with positive coefficients writes
$$X=e_5+2e_7, \ \ Y=2e_6.$$
Thus
$$X \bullet Y=(e_5+2e_7)(4e_6)=4e_5+8e_6=[-12,8].$$
We obtain now the Minkowski product. In general, when one increases the algebra dimension, the product will be closer to the Minkowski one and one still get the distributivity and associativity.

\section{Numerical implementation}
In this section, we show some examples of interval arithmetics applications on simple problems which will prove how this new approach efficient and robust is.

\subsection{Arithmetic implementation in \emph{python}}
We have choosen \emph{python} programming langage\cite{python}. The main reason is that it is a free object-oriented langage, with a huge number of numerical libraries. One of the main advantage of \emph{python} is that first it is possible to link the source code with others written in C/C++, FORTRAN, and second, it interacts easily with other calculations tools such as SAGE\cite{sage} and Maxima\cite{maxima} in order to do formal calculations with \emph{python} langage.
But here, we present pure numerical applications within \emph{python} environnement. The translation in other langages such as C++ and \emph{scilab}\cite{scilab} is very easy and would be available soon.
To start it is necessary to import the $interval\_lib$ library which has been developped to define intervals, vector and matrices of intervals, and all the arithmetic operations.
\begin{verbatim}
Python 2.6.6 (r266:84292, Sep 15 2010, 15:52:39)
[GCC 4.4.5] on linux2
Type "help", "copyright", "credits" or "license" for more information.
>>> from interval_lib import *
\end{verbatim}
The instanciation of an interval $[x,y]$ is done with $interval(x,y,order)$. The variable order corresponds to the dimension of the algebra used to represent the intervals. Its value is set to $4$ by default, which is the minimal one. Another way to define an interval such as $[x-\epsilon,x+\epsilon]$ is $interval(x,eps=\epsilon)$.
Now, let's define the intervals $[-1,2]$, $[3,4]$, $[3,12]$ and $[1,3]$ for example :
\begin{verbatim}
>>> a=interval(-1,2)
>>> b=interval(3,4)
>>> c=interval(3,12)
>>> d=interval(2,eps=1)
\end{verbatim}
It is possible to have more information on each interval :
\begin{verbatim}
>>> print c.min,c.max
3.0 12.0
>>> print abs(c),c.width,c.midpoint
16.5 9.0 7.5
\end{verbatim}
A partial order relation can be implemented on the set of intervals $\mathbb{IR}$ :
\begin{equation}
\forall x,y\in \mathbb{IR}, x\not\subset y,\   center(x)<center(y) \Leftrightarrow x<y
\end{equation}
and
\begin{equation}
\forall x,y\in \mathbb{IR}, x\subset y, width(x)<width(y) \Leftrightarrow x<y
\end{equation}
This can be extended to a total order relation.\\
\begin{verbatim}
>>> print a,b,a<b
[-1.0,2.0] [3.0,4.0] True
\end{verbatim}
\begin{verbatim}
>>> print c,d,d<c
[3.0,12.0] [1.0,3.0] True
\end{verbatim}

\subsection{Semantic and True Arithmetic}
There are two possible arithmetic implementations, depending on the choosen semantic\cite{Jaulin_Kenoufi,Irina_Abdel}. In the first one, called semantic arithmetics, the substraction of two intervals $x$ and $y$ is done according to $x-y=x+(-y)$, and the addition of those two terms. For example $[2,3]-[0,1]=[2,3]+[-1,0]=[2+(-1),3+0]=[1,3]$, and $[-1,1]-[-1,1]=[-1,1]+[-1,1]=[-2,2]\neq[0,0]=0$. As mentionned in the introduction of this paper, "negative" intervals do not have a physical meaning and the addition/substraction between two intervals can not be easily transfered to the bounds of the resulting interval. This yields to the fact that differential calculus in this framework is not relevant and one has to compute the derivatives in the center of the intervals in order to recover a certain meaning. It is not obvious to transfer natural functions to inclusion ones. In the second framework, called true arithmetic, the substractions are done in the algebra $\mathcal{A}_{n}$ with $n\ge 4$.  For the previous example in $\mathcal{A}_{n}$ : $[2,3]-[0,1]=(2,1,0,0)-(0,1,0,0)=(2,0,0,0)=[2,2]$ and for any interval $x$, $x-x=[0,0]=0$. In this arithmetic, it is possible to perform differential calculus and to transfer natural functions to inclusion ones by replacing the terms in the definition by intervals.
One has to note that in both cases, multiplication remains distributive and associative according to addition. But division is distributive for substraction only for the true arithmetic even if it is distributive for addition in the semantic one.
The main reason of this phenomenon, is that the opposite intervals have no real meaning, and it remains to the user to modelize correctly the physical problem. Moreover, there is no wrapping effects and data dependencies as shown on simple examples below.\\
Some other examples :
\begin{verbatim}
>>>  # Semantic arithmetic
>>> print a-a,a*b,b*a,b/b,c+1
[-3.0,3.0] [-4.0,8.0] [-4.0,8.0] [1.0,1.0] [4.0,13.0]
\end{verbatim}
and
\begin{verbatim}
>>>  # True arithmetic
>>> print a-a,a*b,b*a,b/b,c+1
[0.0,0.0] [-4.0,8.0] [-4.0,8.0] [1.0,1.0] [4.0,13.0]
\end{verbatim}
The division is not allowed for intervals containing $0$ :
\begin{verbatim}
>>> print b/a
Interval division in A4 not allowed !!
\end{verbatim}
Let's see the distributive operations :
\begin{verbatim}
>>> # True arithmetic
>>> print a*(b+c),a*b+a*c,(a+b)/c,a/c+b/c
[-16.0,32.0] [-16.0,32.0] [0.5,0.916666666667] [0.5,0.916666666667]

>>> print a*(b-c),a*b-a*c,(a-b)/c,a/c-b/c
[-16.0,8.0] [-16.0,8.0] [-1.08333333333,-0.166666666667] [-1.08333333333,-0.166666666667]
\end{verbatim}
In the semantic arithmetic
\begin{verbatim}
# Semantic arithmetic
>>> print a*(b+c),a*b+a*c,(a+b)/c,a/c+b/c
[-16.0,32.0] [-16.0,32.0] [0.5,0.916666666667] [0.5,0.916666666667]

>>> print a*(b-c),a*b-a*c,(a-b)/c,a/c-b/c
[-28.0,20.0] [-28.0,20.0] [-0.833333333333,-0.416666666667] [-1.08333333333,-0.166666666667]
\end{verbatim}
In the semantic framework, the distributivity of division according substraction is lost but not according addition. This is due to the calculation of $a-b$ before to be divided by $c$. However the division distributivity is always fully respected in the true arithmetic.\\
Another interesting example shows that one gets no wrapping and data dependancy for the two arithmetic frameworks.
\begin{verbatim}
>>> def f1(x):return x**2-2*x+1
>>> def f2(x):return x*(x-2)+1
>>> def f3(x):return (x-1)**2
\end{verbatim}
\begin{verbatim}
>>> # Semantic arithmetic
>>> print f1(a), f2(a), f3(a)
[-7.0,8.0] [-7.0,8.0] [-7.0,8.0]
>>> print f1(b), f2(b), f3(b)
[2.0,11.0] [2.0,11.0] [2.0,11.0]
\end{verbatim}
and
\begin{verbatim}
>>> # True arithmetic
>>> print f1(a), f2(a), f3(a)
[-1.0,2.0] [-1.0,2.0] [-1.0,2.0]
>>> print f1(b), f2(b), f3(b)
[4.0,9.0] [4.0,9.0] [4.0,9.0]
\end{verbatim}
One remarks that the true arithmetic results are always included in the ones obtained with the semantic arithmetic.

\subsection{Optimization examples}

\subsubsection{Fixed-step gradient descent method}
Here is a script example of minimization with fixed-step gradient method which belongs to the so-called gradient descent method\cite{nr}. This algorithm and this example are very simple but it shows that the result is garanted to be found within the final interval.

\begin{figure}[!h]
\centering
\includegraphics[width=6.5cm,height=6.5cm]{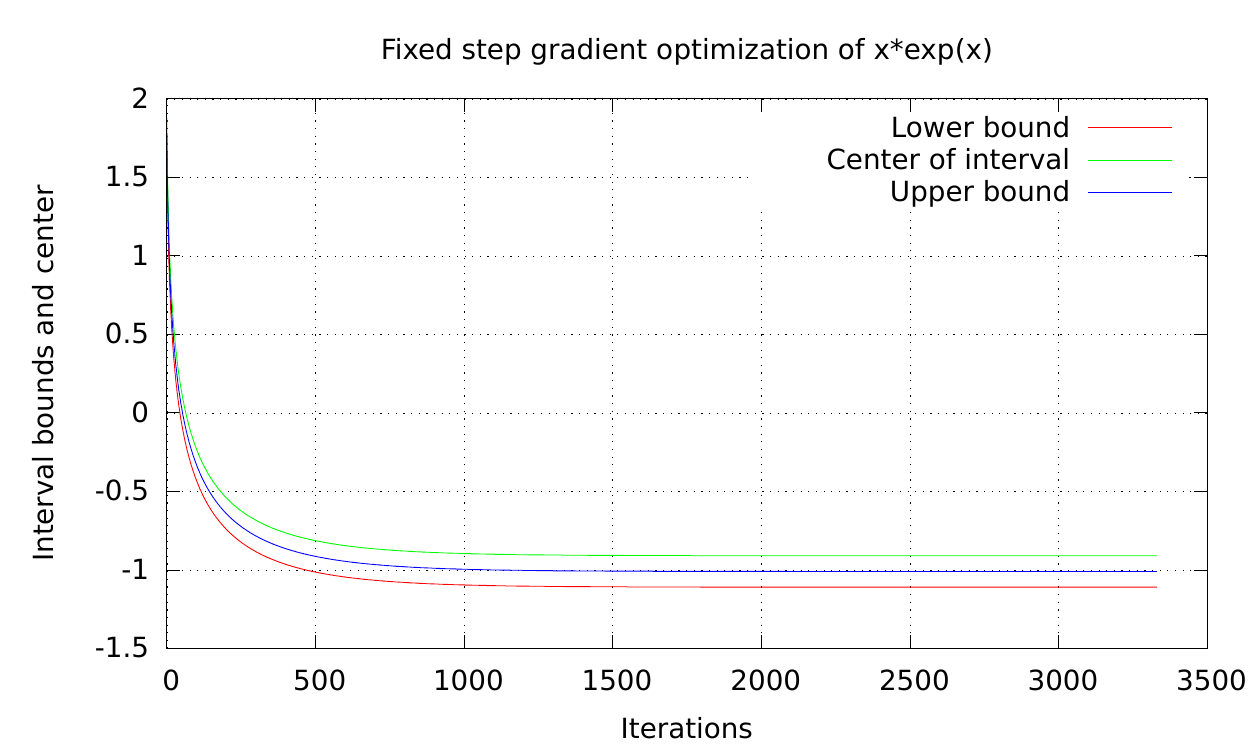}
\caption{Convergence of the fixed-step gradient algorithm for the function $x\mapsto x\cdot\exp(x)$ to an interval centered around $-1$.}
\label{gradient}
\end{figure}

\begin{figure}[!h]
\centering
\includegraphics[width=6.5cm,height=6.5cm]{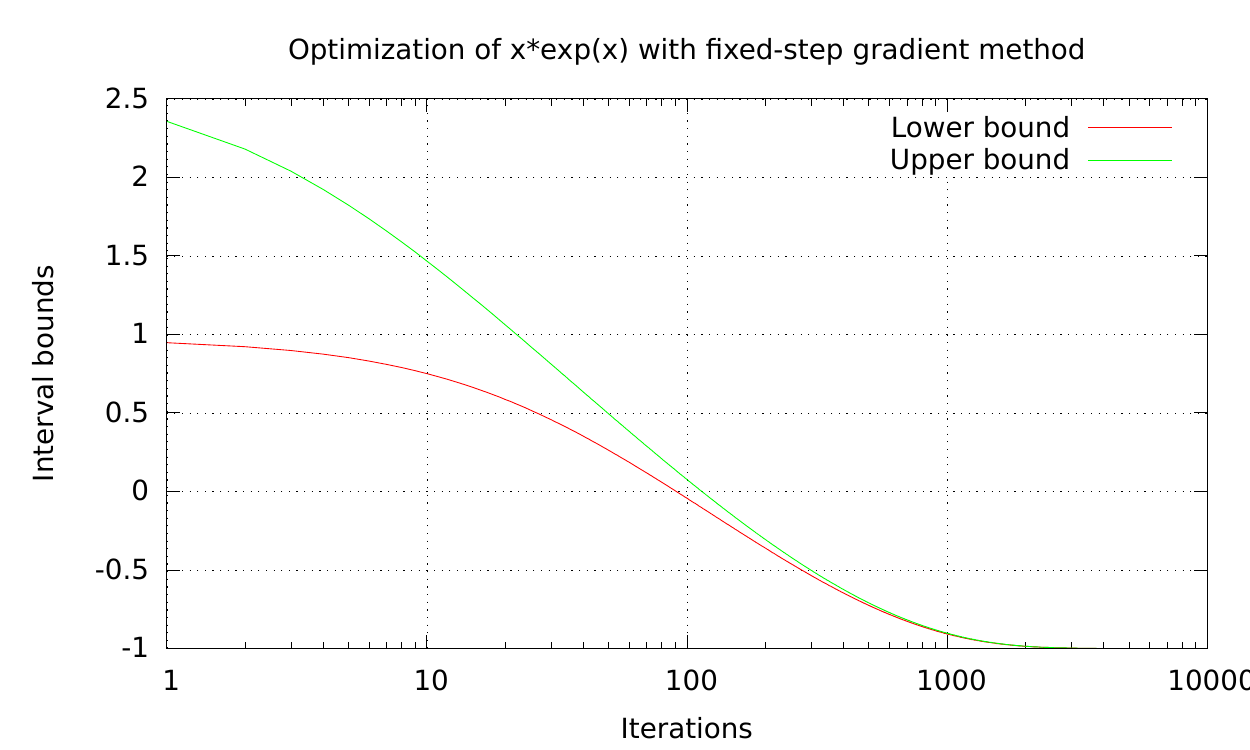}
\caption{Convergence of the fixed-step gradient algorithm with true arithmetic for the function $x\mapsto x\cdot\exp(x)$ to an interval centered around $-1$.}
\label{gradient_true}
\end{figure}

\begin{verbatim}
from interval_lib import *
# Example of fixed step gradient descent method

file=open("res.data", "w") # Data file to be plotted
h=1.e-6 # Finite difference step
def f(x):return x*(exp(x)) # Function to be minimized
def fp(x):return interval(((f(x+h)-f(x-h)).midpoint)/h/2.) # Finite difference
x=interval(2, eps=.1) # Initial guess
rho=interval(1.e-2) # Gradient step
epsilon=1.e-6 # Accuracy of the gradient
while abs(fp(x))>epsilon: # Descent loop
    fprime=fp(x)
    x=x-rho*fprime
    file.write(("%f %f %f %f %f\n")%(x.min, x.max, x.midpoint, fprime.min, fprime.max))
file.close()
\end{verbatim}
In the true arithmetic, the finite differences are "smaller" and it has meaning to do derivative calculations. This is due to the fact that for close intervals, the difference is close to $0$. One has just to change
\begin{verbatim}
def fp(x):return (f(x+h)-f(x-h))/h/2. # Finite difference
\end{verbatim}
The result shown in figure \ref{gradient_true} is impressive, because for any initial guess the interval width decreases to converge to real point minimum. In the semantic interval on figure \ref{gradient}, the width of the interval does not decrease and the center converges to the right value. This is due to finite difference calculation at the center.

\subsubsection{Newton-Raphson method}
Let's optimize the same function $x\mapsto x\cdot\exp(x)$ with a second order method such as the Newton-Raphson one, which is the basis of all second order methods such as Newton or quasi-Newton's ones\cite{nr}. It finds the same minimum which is an interval centered around $-1$.

\begin{verbatim}

# Example of Newton-Raphson method
from interval_lib import *
file=open("res.data", "w") # Data file to be plotted
h=(1.e-6) # Finite difference step
def f(x):return x*(exp(x)) # Function to be minimized
def fp(x):return interval(((f(x+h)-f(x-h)).midpoint)/h/2.) # Finite difference
def fp2(x):return interval(((f(x+h)+f(x-h)-2*f(x)).midpoint)/(h*h)) # Finite difference
x=interval(2, eps=.1) # Initial guess
epsilon=1.e-10 # Accuracy of the gradient
while abs(fp(x))>epsilon: # Descent loop
    fprime=fp(x)
    fsecond=fp2(x)
    file.write(("%f %f %f %f %f\n")%(x.min, x.max, x.midpoint, fprime.min, fprime.max))
    x=x-fprime/fsecond
file.close()
\end{verbatim}

\begin{figure}[!h]
\centering
\includegraphics[width=6.5cm,height=6.5cm]{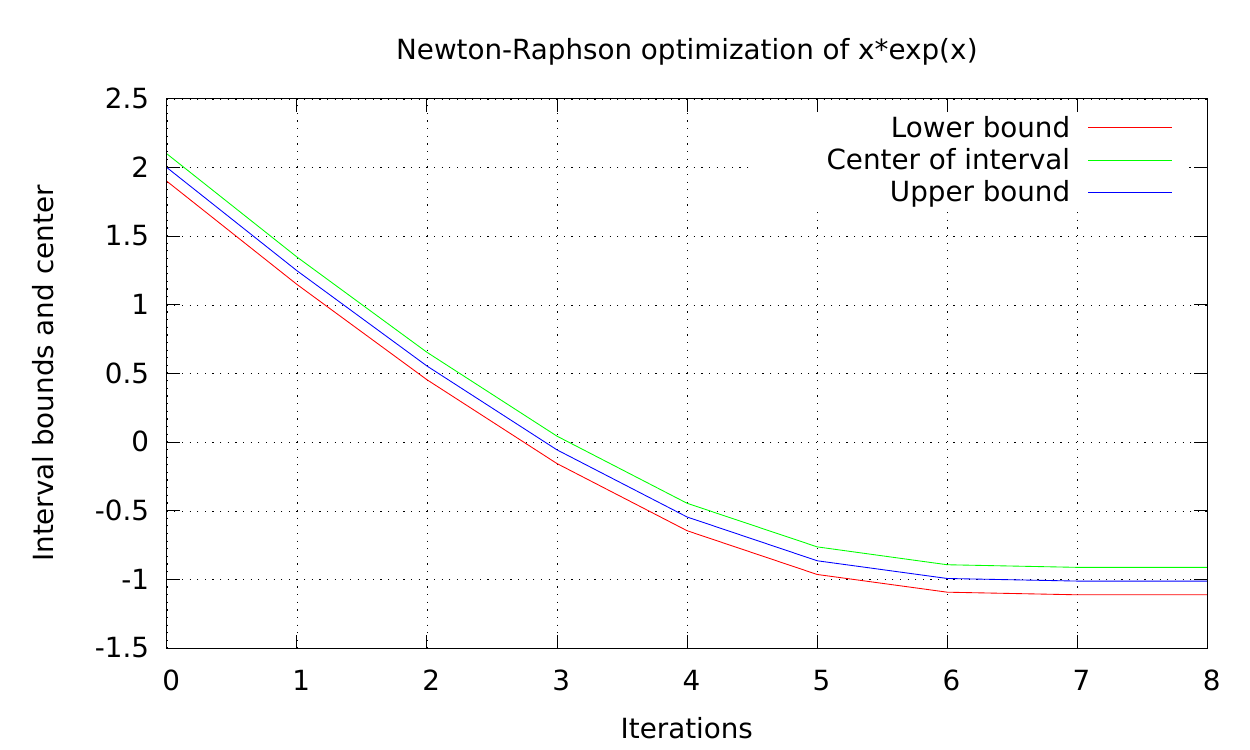}
\caption{Convergence of the Newton-Raphson algorithm for the function $x\mapsto x\cdot\exp(x)$ to an interval centered around $-1$.}
\label{newton}
\end{figure}

\begin{figure}[!h]
\centering
\includegraphics[width=6.5cm,height=6.5cm]{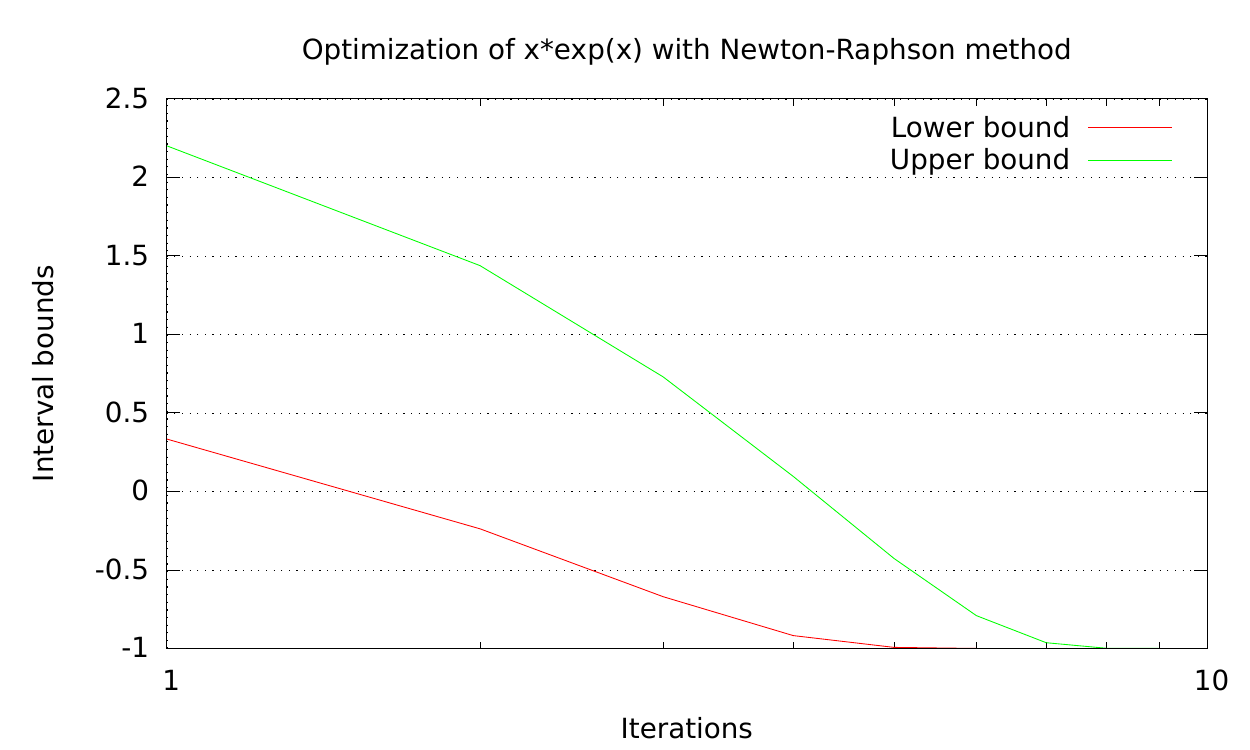}
\caption{Convergence of the Newton-Raphson algorithm with true arithmetic for the function $x\mapsto x\cdot\exp(x)$ to an interval centered around $-1$.}
\label{newton_true}
\end{figure}
\begin{figure}[!h]
\centering
\includegraphics[width=6.5cm,height=6.5cm]{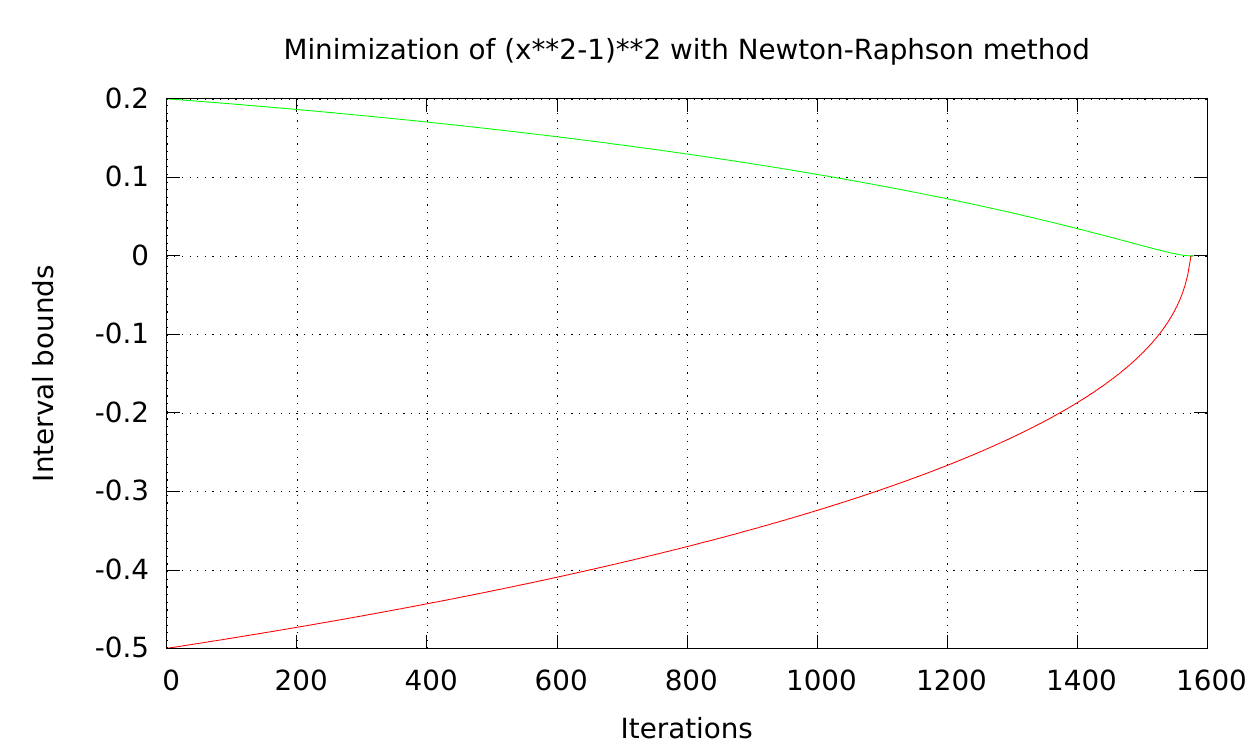}
\caption{Convergence of the Newton-Raphson algorithm with true arithmetic for the function $x\mapsto (x^2-1)^2$.}
\label{newton1}
\end{figure}
\begin{figure}[!h]
\centering
\includegraphics[width=6.5cm,height=6.5cm]{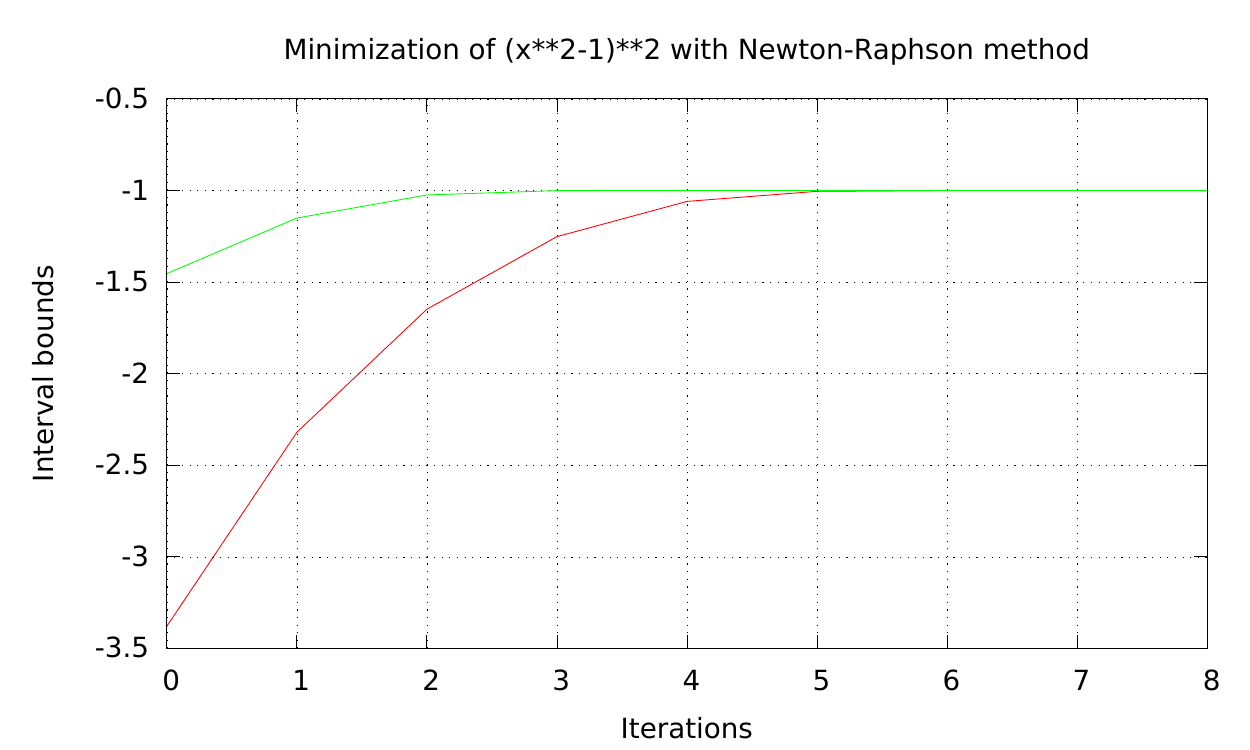}
\caption{Convergence of the Newton-Raphson algorithm with true arithmetic for the function $x\mapsto (x^2-1)^2$.}
\label{newton2}
\end{figure}
\begin{figure}[!h]
\centering
\includegraphics[width=6.5cm,height=6.5cm]{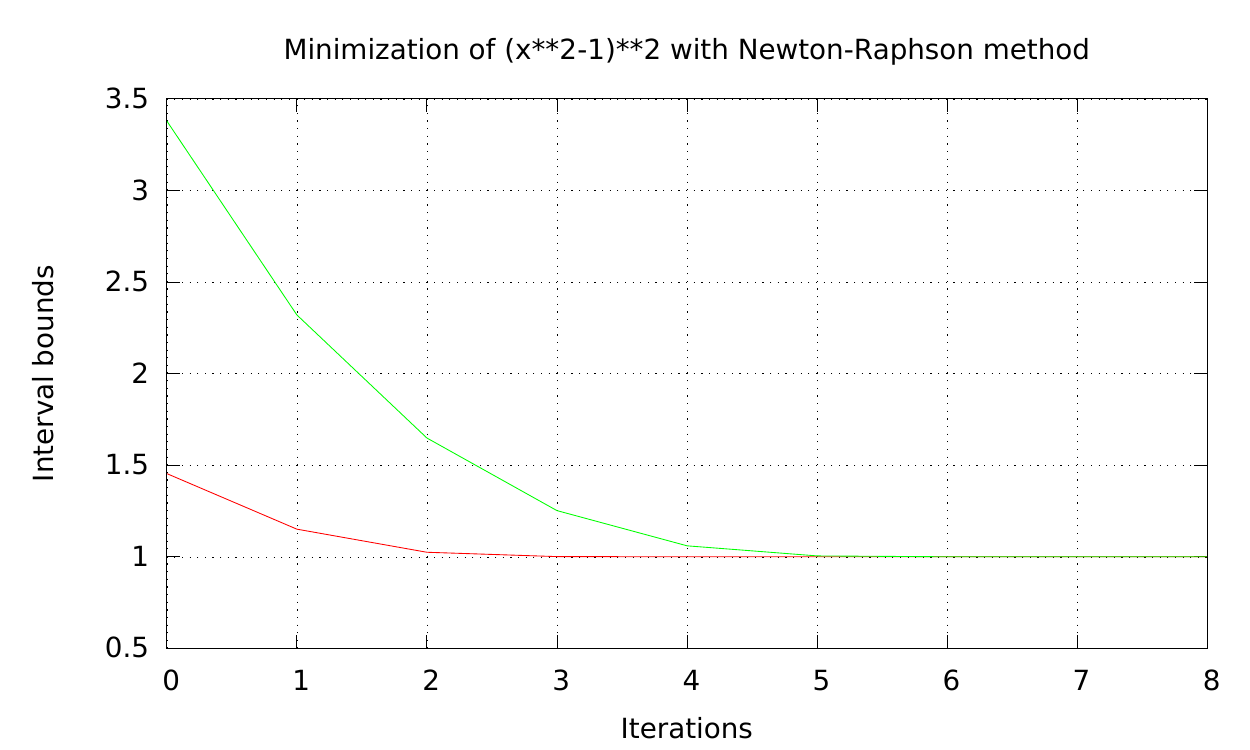}
\caption{Convergence of the Newton-Raphson algorithm with true arithmetic for the function $x\mapsto (x^2-1)^2$.}
\label{newton3}
\end{figure}
Another interesting example is shown on the figures \ref{newton1},\ref{newton2} and \ref{newton3} for different initial guess intervals. We would like to optimize $x\mapsto (x^2-1)^2$. One has to change the finite differences calculated in the center of the interval by classical finite differences :
\begin{verbatim}
def fp(x):return (f(x+h)-f(x-h))/h/2. # Finite difference
def fp2(x):return (f(x+h)+f(x-h)-2*f(x))/(h*h) # Finite difference
\end{verbatim}
The minima are $\{-1,0,1\}$. One can see that depending on the initial guess, this simple algorithm finds the right real point minima.

\begin{figure}[!h]
\centering
\includegraphics[width=6.5cm,height=6.5cm]{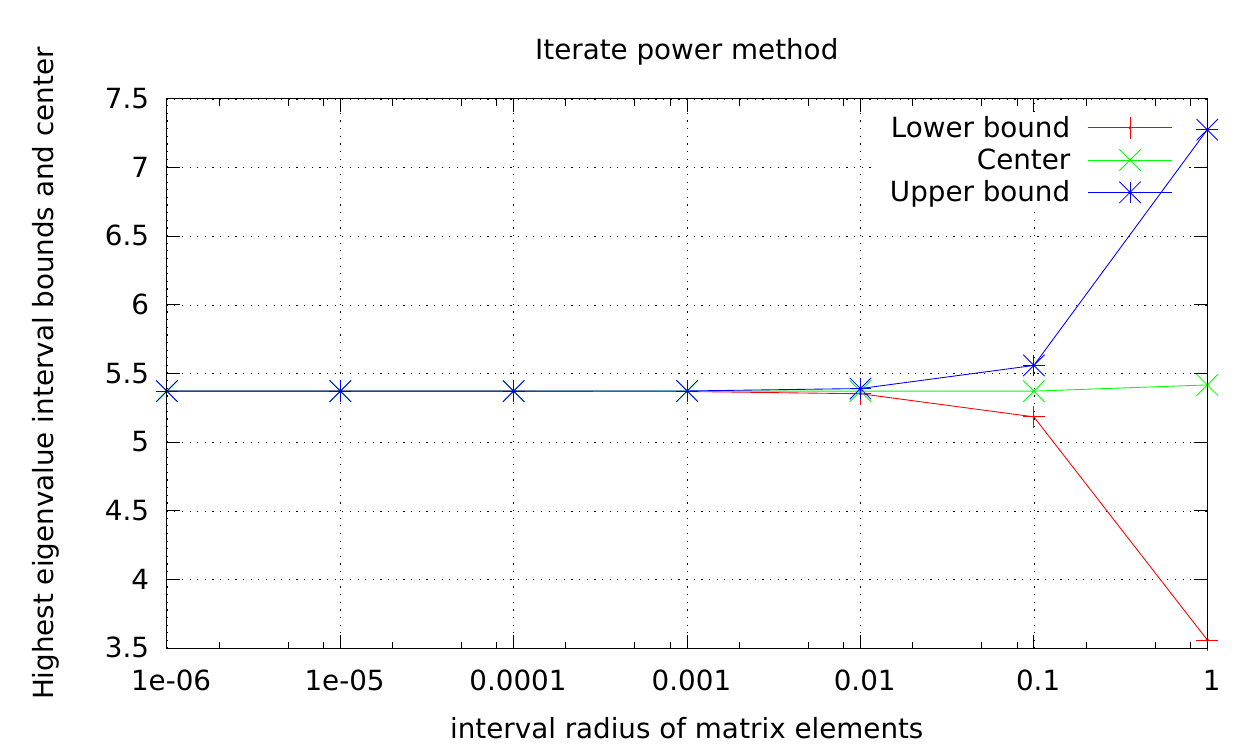}
\caption{Comparison of the largest eigenvalue  found with \emph{scilab} $5.3722813$ versus the iterate power method one computed with intervals.}
\label{Diag1}
\end{figure}

\begin{figure}[!h]
\centering
\includegraphics[width=6.5cm,height=6.5cm]{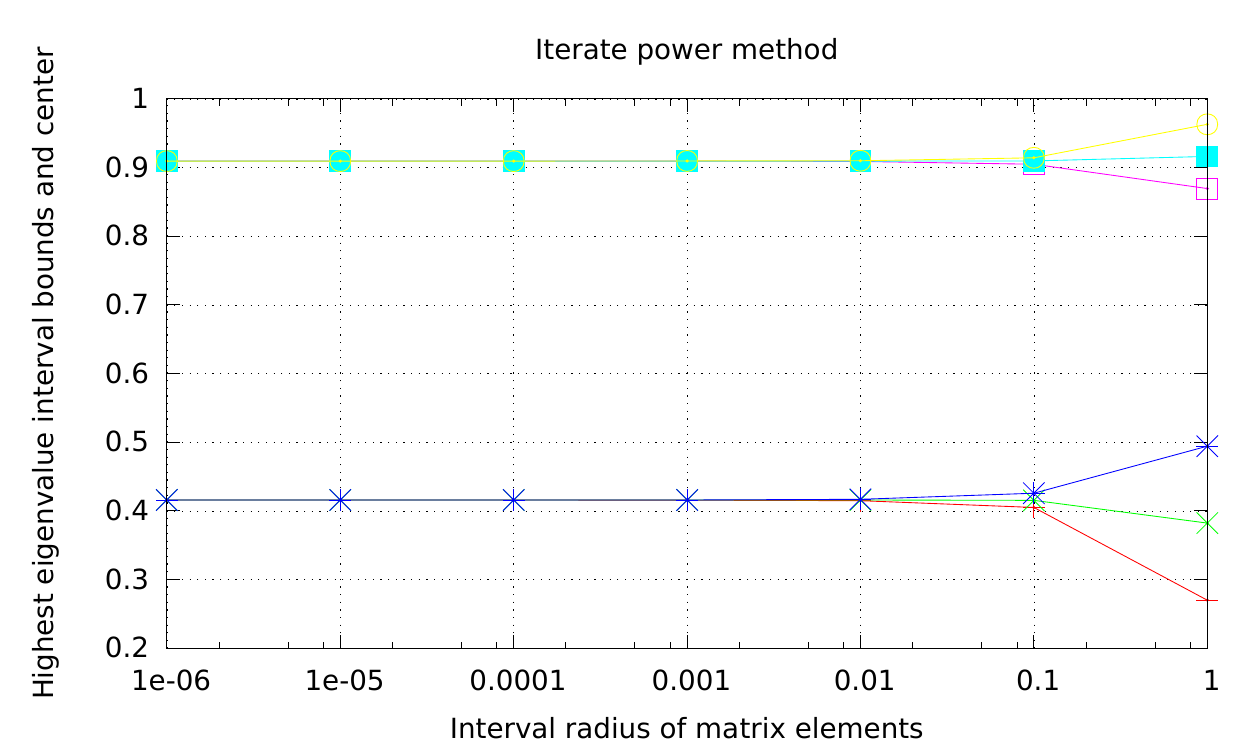}
\caption{Comparison of the eigenvector corresponding to the largest eigenvalue found with \emph{scilab} $(0.4159736,0.9093767)$ versus the iterate power method one computed with intervals.}
\label{Diag2}
\end{figure}

\subsection{Matrix diagonalization and inversion}
\subsubsection{Diagonalization}
As an example, we define the matrix $M$ whose elements are intervals centered around a certain real number with a radius $\epsilon$.
\begin{equation*}
M=\left(
\begin{array}{cc}
\lbrack 1-\epsilon,1+\epsilon] & [2-\epsilon,2+\epsilon] \\
\lbrack 3-\epsilon,3+\epsilon] & [4-\epsilon,4+\epsilon]%
\end{array}%
\right) .
\end{equation*}
If one uses \emph{scilab} to compute the spectrum of the previous matrix without radius ($\epsilon=0$), the highest eigenvalue is approximatively $5.3722813$ and the corresponding eigenvector is $(0.4159736,0.9093767)$.
In order to show that arithmetics and interval algebra developped above is robust and stable, let's try to compute the highest eigenvalue of an interval matrix. One uses here the iterate power method, which is very simple and constitute the basis of several powerful methods such as deflation and others.
The two figures \ref{Diag1} and \ref{Diag2} show clearly for different value of $\epsilon$  the stability of the multiplication, and the largest eigenmode is recovered when $\epsilon=0$. The other eigen modes can be computed with the deflation methods for example which consists to withdraw the direction spanned by the eigenvector associated to the highest eigenvalue to the matrix by constructing its projector and to do the same. Several methods are available and efficient. We have choosen to compute only the highest eigenvalue and its corresponding eigenvector in order to show simply the efficiency of our new artihmetic. The corresponding code in \emph{python} is described below :

\begin{verbatim}
# Example of an interval matrix diagonalization
from interval_lib import *
file=open("res.data", "w") # Data file to be plotted
for i in xrange(10):# Loop on the radius of the matrix elements
    epsilon=10.**(-i)
	# Construction of the matrix
    a=interval(1, eps=epsilon);b=interval(2, eps=epsilon)
    c=interval(3, eps=epsilon);d=interval(4, eps=epsilon)
    u=Vector([a, b]);v=Vector([c, d])
    u0=Vector([interval(1), interval(1)]) #Initial guess
    m=Matrix([u, v]) # Interval matrix to be diagonalized
    e,v=iterate_power(m, u0, 10) #Power iteration
    file.write("%f %f %f %f %f %f %f\n"%(epsilon, e.min, e.max, v[0].min, v[0].max,v[1].min, v[1].max ))
file.close()
\end{verbatim}

\subsubsection{Inversion}
Let's define a symmetric matrix.
\begin{equation*}
M=\left(
\begin{array}{ccc}
\lbrack 1-\epsilon,1+\epsilon] & [4-\epsilon ,4+\epsilon] & [5-\epsilon, 5+\epsilon] \\
\lbrack 4-\epsilon, 4+\epsilon] & [2-\epsilon,2+\epsilon] & [6-\epsilon, 6+\epsilon]\\
\lbrack 5-\epsilon, 5+\epsilon] & [5-\epsilon, 6+\epsilon]& [3-\epsilon, 3+\epsilon]
\end{array}%
\right) .
\end{equation*}
We would like to use the well-know Schutz-Hotelling algorithm\cite{nr,Householder} to inverse a matrix $X$ :
\begin{equation}
X_0=\frac{X^t}{\sum_{i,j}A_{ij}^2},\ X_j=X_{j-1}(2-A\cdot X_{j-1}),\ \forall n\ge 1
\end{equation}
\emph{Scilab} gives numerically for $\epsilon=0$
\begin{equation}
M^{-1}=\left(
\begin{array}{ccc}
  - 0.2678571   & 0.1607143   & 0.125  \\
    0.1607143 & - 0.1964286 &   0.125  \\
    0.125     &   0.125  &    - 0.125
    \end{array}%
\right).
\end{equation}
The \emph{python} code is very simple :
\begin{verbatim}
# Example of an interval matrix inversion
from interval_lib import *

epsilon=.2 #intervals radius
a=interval(1, eps=epsilon);b=interval(2, eps=epsilon);c=interval(3, eps=epsilon)
d=interval(4, eps=epsilon);e=interval(5, eps=epsilon);f=interval(6, eps=epsilon)
# Build the matrix m
u=Vector([a, d, e]);v=Vector([d, b, f]);w=Vector([e, f, c])
m=Matrix([u, v, w]);inv_m=schultz(m);inv_inv_m=schultz(inv_m)
# Display results
print "M=", m
print "Inverse matrix = ", inv_m
print "M^(-1)*M=", inv_m*m
print "M*M^(-1)=", m*inv_m
print "(M^(-1))^(-1)=",inv_inv_m
\end{verbatim}
We obtain with intervals for $\epsilon=0.2$ for example :
\begin{verbatim}
M= [*
[0.8,1.2][3.8,4.2][4.8,5.2]
[3.8,4.2][1.8,2.2][5.8,6.2]
[4.8,5.2][5.8,6.2][2.8,3.2]
*]
Inverse matrix =  [*
[-0.267918088737,-0.267790262172][0.160409556314,0.161048689139][0.12457337884,0.125468164794]
[0.160409556314,0.161048689139][-0.19795221843,-0.194756554307][0.122866894198,0.12734082397]
[0.12457337884,0.125468164794][0.122866894198,0.12734082397][-0.127986348123,-0.121722846442]
*]
M^(-1)*M= [*
[1.0,1.0][-2.77555756156e-17,0.0][-1.80411241502e-16,-1.66533453694e-16]
[-2.22044604925e-16,-2.11636264069e-16][1.0,1.0][-1.38777878078e-17,0.0]
[0.0,1.38777878078e-17][0.0,0.0][1.0,1.0]
*]
M*M^(-1)= [*
[1.0,1.0][-1.21430643318e-16,-1.11022302463e-16][0.0,1.38777878078e-17]
[0.0,0.0][1.0,1.0][-1.38777878078e-17,0.0]
[4.16333634234e-17,5.55111512313e-17][-1.80411241502e-16,-1.66533453694e-16][1.0,1.0]
*]
(M^(-1))^(-1)= [*
[0.8,1.2][3.8,4.2][4.8,5.2]
[3.8,4.2][1.8,2.2][5.8,6.2]
[4.8,5.2][5.8,6.2][2.8,3.2]
*]
\end{verbatim}
and for $\epsilon=0.1$
\begin{verbatim}
M= [*
[0.9,1.1][3.9,4.1][4.9,5.1]
[3.9,4.1][1.9,2.1][5.9,6.1]
[4.9,5.1][5.9,6.1][2.9,3.1]
*]
Inverse matrix =  [*
[-0.267888307155,-0.267824497258][0.160558464223,0.160877513711][0.124781849913,0.125228519196]
[0.160558464223,0.160877513711][-0.197207678883,-0.195612431444][0.123909249564,0.126142595978]
[0.124781849913,0.125228519196][0.123909249564,0.126142595978][-0.126527050611,-0.123400365631]
*]
M^(-1)*M= [*
[1.0,1.0][2.22044604925e-16,2.35922392733e-16][1.66533453694e-16,1.7694179455e-16]
[0.0,0.0][1.0,1.0][0.0,0.0]
[-1.31838984174e-16,-1.11022302463e-16][-1.2490009027e-16,-1.11022302463e-16][1.0,1.0]
*]
M*M^(-1)= [*
[1.0,1.0][0.0,0.0][-6.93889390391e-18,0.0]
[0.0,0.0][1.0,1.0][0.0,6.93889390391e-18]
[-3.46944695195e-18,0.0][0.0,0.0][1.0,1.0]
*]
(M^(-1))^(-1)= [*
[0.9,1.1][3.9,4.1][4.9,5.1]
[3.9,4.1][1.9,2.1][5.9,6.1]
[4.9,5.1][5.9,6.1][2.9,3.1]
*]
\end{verbatim}
and for $\epsilon=0.01$
\begin{verbatim}
M= [*
[0.99,1.01][3.99,4.01][4.99,5.01]
[3.99,4.01][1.99,2.01][5.99,6.01]
[4.99,5.01][5.99,6.01][2.99,3.01]
*]
Inverse matrix =  [*
[-0.267860324247,-0.267853946662][0.160698378764,0.160730266691][0.124977730269,0.125022373367]
[0.160698378764,0.160730266691][-0.196508106182,-0.196348666547][0.124888651345,0.125111866834]
[0.124977730269,0.125022373367][0.124888651345,0.125111866834][-0.125155888117,-0.124843386433]
*]
M^(-1)*M= [*
[1.0,1.0][1.11022302463e-16,1.12323345069e-16][5.55111512313e-17,5.63785129692e-17]
[-1.51788304148e-18,0.0][1.0,1.0][-5.74627151417e-17,-5.55111512313e-17]
[-8.67361737988e-19,0.0][-2.22911966663e-16,-2.22044604925e-16][1.0,1.0]
*]
M*M^(-1)= [*
[1.0,1.0][-2.16840434497e-19,0.0][0.0,0.0]
[-2.22044604925e-16,-2.21610924056e-16][1.0,1.0][0.0,0.0]
[-1.12323345069e-16,-1.11022302463e-16][-5.57279916658e-17,-5.55111512313e-17][1.0,1.0]
*]
(M^(-1))^(-1)= [*
[0.99,1.01][3.99,4.01][4.99,5.01]
[3.99,4.01][1.99,2.01][5.99,6.01]
[4.99,5.01][5.99,6.01][2.99,3.01]
*]
\end{verbatim}

It is obvious that this method is very stable and confirms that the true arithmetic operations are robust. It is not difficult to extend usual linear iterative algebra numerical algorithms to intervals. It permits to solve a lot problems where the entries of the matrices are not well defined, especially for automation applications\cite{Irina_Abdel}.

\section{Conclusion}
We have presented a better algebraic way to do calculations on intervals. This approach \cite{nicolas} is done by embedding the space of intervals into a free algebra of dimension greater or equal to $4$. This permits to obtain all the basic arithmetic operators with distributivity and associativity. We have shown that when one increases the representative algebra dimension, the multiplication result will be closer to the usual Minkowski product. We have compared two approaches for interpreting substraction operation, and the canonical approach we have proposed, called true arithmetics is more coherent and efficient. Differential calculus is possible and very efficient to solve some optimization problems. It is now possible to build inclusion functions from the natural ones. This will be studied in a more accurate way in a forthcoming paper. The set of intervals is now endowed with an order relation which permits to define inequalities for intervals. One of the straightforward application can be non-linear simplex algorithm, the so-called Nelder-Mead simplex method or downhill simplex\cite{nr,nm} which is derivative-free method and can be easily implemented.
We have exhibited some examples of applications : optimization, diagonalization and inversion of matrices which clearly state that the arithmetic is stable and that if the initial datas are known with a certain uncertainity (belonging to an interval), it is thus possible to estimate with accuracy the point solution of the problem, a real number or a smaller interval centered around it.
\begin{acknowledgements}
We thank Michel Gondran and Irina Berseneva for useful and interesting discussions.
\end{acknowledgements}

%\tableofcontents

\end{document}